
\input amstex
\documentstyle{amsppt}
\loadbold
\magnification=\magstep1
\vsize 7in
\baselineskip18pt

\hyphenation{co-deter-min-ant co-deter-min-ants pa-ra-met-rised
pre-print pro-pa-gat-ing pro-pa-gate
fel-low-ship Cox-et-er dis-trib-ut-ive}
\def\leaderfill{\leaders\hbox to 1em{\hss.\hss}\hfill}

\def\a{{\alpha}}

\def\e{{\varepsilon}}

\def\l{{\lambda}}

\def\bu{{\bold u}}

\def\bw{{\boldsymbol \omega}}
\def\bx{{\bold x}}

\def\b0{\text{\bf 0}}

\def\ra{{\ \longrightarrow \ }}

\def\supp{\text{\rm \, supp}}

\def\zed{{\Bbb Z}}

\def\enn{{\Bbb N}}

\def\End{\text{\rm End}}
\def\Hom{\text{\rm Hom}}

\def\boxit#1{\vbox{\hrule\hbox{\vrule \kern3pt
\vbox{\kern3pt\hbox{#1}\kern3pt}\kern3pt\vrule}\hrule}}
\def\rabbit{\vbox{\hbox{\kern0pt
\vbox{\kern0pt{\hbox{---}}\kern3.5pt}}}}

\def\tableau#1{
        \hbox {
                \hskip -10pt plus0pt minus0pt
                \raise\baselineskip\hbox{
                \offinterlineskip
                \hbox{#1}}
                \hskip0.25em
        }
}

\def\tabCol#1{
\hbox{\vtop{\hrule
\halign{\strut\vrule\hskip0.5em##\hskip0.5em\hfill\vrule\cr\lower0pt
\hbox\bgroup$#1$\egroup \cr}
\hrule
} } \hskip -10.5pt plus0pt minus0pt}

\def\CR{
        $\egroup\cr
        \noalign{\hrule}
        \lower0pt\hbox\bgroup$
}



\def\blank#1#2{
\hbox to #1{\hfill \vbox to #2{\vfill}}
}


\def\strut{\vrule height10pt depth5pt width0pt}

\def\seca{1}
\def\secb{2}
\def\secba{3}
\def\secc{4}
\def\secd{5}
\def\sece{6}
\def\secf{7}
\def\asn{T(n)}
\def\asnt{T(n)[Q^{-1}]}
\def\asm{M_q(n)}
\def\asmt{\widetilde{M}_q(n)}
\def\asml#1#2{M_{ {#1}, {#2} }(n)}
\def\asmlt#1#2{\widetilde{M}_{ {#1}, {#2} }(n)}
\def\kq{{k[q]}}
\def\kqq{{k[q,q^{-1}]}}
\def\hcat{\bold{Heap}}
\def\hcatg{{\bold{Heap}(\Gamma)}}
\def\bw{\bold{w}}
\def\bx{\bold{x}}
\def\bu{\bold{u}}
\def\heapel{{E(n)}}
\def\cog{\text{Co}(\Gamma)}
\def\Mod{\text{Mod}}
\def\supp{{\text{supp}}}
\def\propid{{{\Cal B}(E(n))}}.
\def\twd{{D_{n+1}^{(2)}}}
\def\mlq{{M_\Lambda[q]}}
\def\mlqq{{M_\Lambda[q, q^{-1}]}}

\topmatter
\title The nil Temperley--Lieb algebra of type affine $C$
\endtitle

\author R.M. Green \endauthor
\affil Department of Mathematics \\ University of Colorado Boulder \\
Campus Box 395 \\ Boulder, CO  80309-0395 \\ USA \\ {\it  E-mail:}
rmg\@euclid.colorado.edu \\
\newline
\endaffil


\abstract
We introduce a type affine $C$ analogue of the
nil Temperley--Lieb algebra,
in terms of generators and relations. We show that this algebra $\asn$, which is
a quotient of the positive part of a Kac--Moody algebra of type
$\twd$, has an easily described faithful representation as an algebra
of creation and annihilation operators on particle configurations, reminiscent
of the open TASEP model in statistical physics. The centre of $\asn$
consists of polynomials in a certain element $Q$, and $\asn$ is a free module
of finite rank over its centre. We show how to localize $\asn$ by adjoining
an inverse of $Q$, and prove that the resulting algebra is a full matrix 
ring over a ring of Laurent polynomials over a field. Although $\asn$ has
wild representation type, over an algebraically closed field we can 
classify all the finite dimensional indecomposable representations of
$\asn$ in which $Q$ acts invertibly.
\endabstract


\endtopmatter


\head \seca. Introduction \endhead

The nil Temperley--Lieb algebra of type affine $A$ 
is an infinite dimensional associative
algebra given by generators and relations. The algebra was introduced by
Postnikov \cite{{\bf 15}} who used it to study the quantum cohomology of 
Grassmannians, but it also has connections to
Weyl algebras, Clifford algebras, and universal enveloping algebras.
(More details on these connections may be found in the introduction of
\cite{{\bf 1}} and references therein.) The nil Temperley--Lieb algebra 
of type affine $A$ can also be constructed as an algebra
of operators on fermionic particle configurations on a circle (described
in \cite{{\bf 13}}), or as the
associated graded algebra of the affine Temperley--Lieb algebra.

Nil Temperley--Lieb algebras have been defined for Coxeter systems of
all types by Biagioli, Jouhet and Nadeau in \cite{{\bf 2}}. In each case, the 
algebra has a
basis indexed by the fully commutative elements of the Weyl group, which were
introduced by Stembridge \cite{{\bf 16}}. Biagioli et al point out that these 
algebras have not been studied much, except in the cases of Weyl groups
of type affine $A$ (as discussed above) and the case of the symmetric group.
The algebra in the symmetric group case is known simply as the 
nil Temperley--Lieb algebra, and it was originally introduced by
Fomin and Greene \cite{{\bf 6}} in their study of symmetric functions, motivated
in turn by work of Billey, Jockusch and Stanley \cite{{\bf 3}}.

In this paper we define, using generators and relations, a 
nil Tem\-per\-ley--Lieb algebra of type affine $C$, or ``nTL algebra
of type affine $C$'' for short. This
algebra has a basis indexed by the minuscule elements of the Weyl group 
rather than the fully commutative elements. In types $A$ and affine $A$, the
fully commutative elements and minuscule elements coincide, but in general
they do not, which means that our algebra is a proper quotient of the 
generalized nil Temperley--Lieb algebra introduced by Biagioli et al.
We study this smaller algebra because it turns out to have much more favourable 
properties from the point of view of representation theory and combinatorics.

In particular, we will show how the nTL algebra of type affine $C$ can
be faithfully constructed in terms of an algebra of creation and
annihilation operators on suitably defined particle configurations.
A large part of the appeal of the nTL algebra
$\asn$ comes from the fact that it can be described so concisely as an
algebra of operators. (This point of view can also 
be used to realize the algebra $\asn$ as a subalgebra of a
$q$-deformation of a certain Clifford algebra, although we do not do this
here.) These particle configurations, and the operators that act on them,
also appear in
the open Totally Asymmetric Simple Exclusion Process (TASEP), which arises in
statistical physics. The open TASEP is also used in molecular biology, where
it is used to model the process of translation in protein synthesis.
For further details of the connections to molecular
biology and statistical physics, the reader is referred to the survey paper
\cite{{\bf 19}} and the references therein.

The proof that the algebra of operators gives a faithful representation of
the nTL algebra of type affine $C$ relies on some combinatorial properties
of reduced words in the corresponding Weyl group.
Although the particle configuration representation is very easy
to describe (see Proposition \secb.3), it is difficult to work with in proofs.
It turns out to be much easier for our purposes
to use the framework of heaps, which are certain
labelled partially ordered sets.

The key heap-theoretic property (proved in Theorem \secc.13) is that every
minuscule element can be associated to a convex subheap of a certain heap;
as we will explain in our concluding remarks, this corresponds to a certain 
property of Coxeter elements.
As a by-product, we prove a general result (Theorem \secc.1) that
shows how under certain hypotheses, convex subheaps can be characterized
by local information; this result seems to be new.

Using the heap approach, we prove that the centre of $\asn$ over a field $k$
consists of polynomials in a certain element $Q$, and that 
$\asn$ is a free module of finite rank over its centre. Adjoining an inverse
of $Q$ to $\asn$ produces the full ring of $2^n$ by $2^n$ matrices over
the ring $\kqq$ of Laurent polynomials. If $k$ is algebraically closed,
the Morita equivalence of this ring
with $\kqq$ can be exploited to classify 
all the finite dimensional indecomposable representations of
$\asn$ in which $Q$ acts invertibly. 
The other finite dimensional indecomposable representations are those in 
which $Q$ acts nilpotently, but classifying these representations is provably
hopeless, even in the case where $Q$ acts as zero,
because the quotient algebra $\asn/\langle Q \rangle$ has wild representation
type (see Theorem \secd.11 (iii)).

If we identify the generators $u_i$ with the 
corresponding generators $E_i$ of the positive part of the affine Kac--Moody
algebra ${\frak g}$ 
of type $\twd$, comparison of our presentation with the Serre 
presentation shows that $\asn$ is a quotient of the plus part of the universal
enveloping algebra $U({\frak g})$ of ${\frak g}$. Although we do not pursue
this here, the full matrix ring of the previous paragraph can be made into 
a representation of the loop algebra of ${\frak g}$ in a similar way.

The nil Temperley--Lieb algebra of type
affine $C$ could also be constructed as the graded algebra of a 
version of the symplectic blob algebra. The latter algebra was
introduced by Martin, Parker and the author \cite{{\bf 14}}, and independently by 
de Gier and Nichols \cite{{\bf 7}}, and has applications to statistical mechanics.
In an analogous way, the type affine $C$ algebra of Biagioli et al can
be regarded as the graded version of the generalized Temperley--Lieb algebra
of type affine $C$ in the sense of Graham \cite{{\bf 8}}; the latter algebra
has been studied in detail by Ernst \cite{{\bf 4}, {\bf 5}}.

The construction of the algebra of operators on particle configurations can
be performed more generally using a minuscule representation for a simple Lie
algebra as a starting point, and in particular, the fermionic particle 
representations of
\cite{{\bf 13}} can be constructed from the minuscule representations in type 
$A_{n-1}$. The representation governing the combinatorics in the 
nTL algebra of type affine $C$ is the spin representation in type $B_n$. 
This is a representation of dimension $2^n$, which is (not coincidentally) the
dimension of most of the irreducible modules for the nTL algebra of type 
affine $C$. We hope to explore these more general constructions in future work.

\head \secb. Definitions \endhead

In Section \secb, we will recall some important definitions and state most
of our main results.

A {\it Coxeter system} is a pair $(W, S)$, where $W$ is a group given by
the presentation $$
\langle S \ | \ (s_i s_j)^{m_{ij}} = 1 \rangle
,$$ where $m_{ij} \in \zed \cup \{\infty\}$ and we have $m_{ii} = 1$ and
$m_{ij} = m_{ji}$. If $m_{ij} = \infty$, we omit the corresponding relation.

Biagioli, Jouhet and Nadeau \cite{{\bf 2}, \S6} associate to an arbitrary
Coxeter system $(W, S)$
a unital associative algebra, nTL$(W)$, 
called the {\it generalized nil Temperley--Lieb 
algebra}. It is given by generators $\{1\} \cup \{u_s: s \in S\}$ and defining 
relations $$\eqalign{
u_s^2 &= 0, \cr
u_s u_t &= u_t u_s \text{\quad if } m(s, t) = 2, \cr
\underbrace{u_s u_t u_s u_t \cdots}_{m_{s, t}} &= 0 \text{\quad if } 2 < m(s, t) 
< \infty. \cr
}$$ 

The Coxeter systems in this paper will all be of type affine $C_n$ unless 
otherwise stated. This means
that we can take $S = \{0, 1, \ldots, n\}$, and $$
m_{s, t} = \cases
1 & \quad \text{ if } s = t,\cr
2 & \quad \text{ if } |s - t| > 1,\cr
3 & 
\quad \text{ if } |s - t| = 1 \text{ and } 1 \leq s, t \leq n-1, 
\text{ and }\cr
4 & \quad \text{ if } \{s, t\} = \{0, 1\} \text{ or } \{s, t\} = \{n-1, n\}. \cr
\endcases
$$ In this special case, the algebra nTL$(W)$ is given by generators
$\{1\} \cup \{u_i : 0 \leq i \leq n\}$ and relations $$\eqalign{
u_i^2 &= 0 \quad \text{ for all } i,\cr
u_i u_j &= u_j u_i \quad \text{ if } |i - j| > 1,\cr
u_i u_j u_i &= 0,
\quad \text{ if } |i - j| = 1 \text{ and } 1 \leq i, j \leq n-1,\cr 
u_0 u_1 u_0 u_1 &= u_1 u_0 u_1 u_0 = 0, \text{ and}\cr
u_{n-1} u_n u_{n-1} u_n &= u_n u_{n-1} u_n u_{n-1} = 0.\cr
}$$

The type affine $C$ nil Temperley--Lieb algebra is
obtained by modifying the last two relations above, 
as follows.

\definition{Definition \secb.1}
The type affine $C$ nil Temperley--Lieb algebra, $\asn$, over a field $k$
is the associative
unital $k$-algebra given by generators $\{1\} \cup \{u_i : 0 \leq i \leq n\}$ 
and defining relations $$\eqalign{
u_i^2 &= 0 \quad \text{ for all } i,\cr
u_i u_j &= u_j u_i \quad \text{ if } |i - j| > 1,\cr
u_i u_j u_i &= 0,
\quad \text{ if } |i - j| = 1 \text{ and } 1 \leq i, j \leq n-1,\cr 
u_1 u_0 u_1 &= 0, \text{ and}\cr
u_{n-1} u_n u_{n-1} &= 0.\cr
}$$
\enddefinition

It follows immediately from the definitions that $\asn$ is a quotient of
nTL$(W)$ if $W$ has type affine $C$.

A key tool for understanding the structure of $\asn$
is a module, $\asm$, on which it acts.

\definition{Definition \secb.2}
Let $k$ be a field and let $q$ be an indeterminate. Let $\asm$ be the free
$\kq$-module with basis consisting of all $2^n$ length $n$ strings on the 
alphabet $\{+, -\}$, and define $\asmt := \kqq \otimes_\kq \asm$.
\enddefinition

\proclaim{Proposition \secb.3}
The $k$-vector space $\asm$ becomes a left $\asn$-module in which 
$\asn$ acts via the following 
$\kq$-linear transformations, where $i$ satisfies 
$1 \leq i \leq n-1$:

\vskip 10pt \quad\quad\quad\quad $\displaystyle{u_0(r_0 r_1 \cdots r_{n-1}) =
\cases
q.({-}r_1 \cdots r_{n-1}) & \text{ if } r_0 = {+},\cr
0 & \text{ otherwise;}\cr
\endcases}$

\vskip 10pt \quad\quad\quad\quad $\displaystyle{u_i(r_0 r_1 \cdots r_{n-1}) = 
\cases
r_0r_1 \cdots r_{i-2}{+}{-}r_{i+1} \cdots r_{n-1} & \text{ if } 
r_{i-1}r_i = {-}{+},\cr
0 & \text{ otherwise;}\cr
\endcases}$

\vskip 10pt \quad\quad\quad\quad $\displaystyle{u_n(r_0 r_1 \cdots r_{n-1}) = 
\cases
r_0 \cdots r_{n-2}{+} & \text{ if } r_{n-1} = {-},\cr
0 & \text{ otherwise.}\cr
\endcases}$

A similar result holds for $\asmt$, after extending scalars to $\kqq$.
\endproclaim

\demo{Proof}
This is a matter of checking that the defining relations in Definition
\secb.1 are satisfied. We omit the proof of this for reasons of space.
\qed\enddemo

Proposition \secb.3 can be thought of in terms of creation and annihilation
operators of particle configurations
if we interpret the symbols ${+}$ (respectively, ${-}$) as the presence 
(respectively, absence) of a particle in a particular position. From this
point of view, $u_0$ acts as a scaled annihilation operator and $u_n$ acts 
as a creation operator. The other generators $u_i$ act by a combination
of annihilating the particle in position $i$ and creating a particle in 
position $i-1$; alternatively, we can think of this as moving a particle in
position $i$ (if there is one) to an empty spot in position $i-1$ 
(if there is one). These operators on particle configurations also arise
in the open TASEP model in statistical physics and molecular
biology, although in the scientific literature (as described
in \cite{{\bf 19}}) the particles
are typically considered to move from left to right.

Two of our main results will be the following.

\proclaim{Theorem \secb.4}
The $k$-algebra $\asn$ acts faithfully on $\asm$, and on $\asmt$.
\endproclaim

Theorem \secb.4, which will be proved at the end of Section \secd, 
will mean that we can regard $\asn$ as an algebra of 
operators on $\asm$, or on $\asmt$. 

With this identification, we have the following result, which will be proved
in Section \secf.

\proclaim{Theorem \secb.5}
There exists an element $Q \in \asn$ whose action on $\asm$ is given by
multiplication by $q$. The centre of $\asn$ is equal to $k[Q]$, and 
$\asn$ is a free $k[Q]$-module of rank $2^{2n}$.
\endproclaim

Since the action of $Q$ on $\asmt$ is invertible, we can extend the algebra
$\asn$ by adjoining a (central) two-sided inverse of $Q$. The next result,
which is proved in Section \sece, shows that the resulting
algebra $\asnt$ turns out to have a remarkably simple structure.

\proclaim{Theorem \secb.6}
The $\kqq$-algebra $\asnt$, in which $q$ acts by multiplication by $Q$,
is isomorphic to the full matrix ring
$M_{2^n}(\kqq)$.
\endproclaim

If $k$ is algebraically closed and $M$ is a finite dimensional $\asn$-module,
it follows by considering the Jordan canonical form of the action of $Q$ 
on $M$ that $M$ is the direct sum of the generalized eigenspaces for the
action of $Q$. Since $Q$ is central, these generalized
eigenspaces are $\asn$-submodules.
In particular, if $M$ is indecomposable, the action of $Q$ on $M$ is either
invertible or nilpotent.

Theorem \secb.6 implies that $\asnt$ is Morita equivalent to the
principal ideal domain $\kqq$. This makes it easy to classify the finite
dimensional indecomposable modules for $\asnt$, as well as the finite 
dimensional
indecomposable modules for $\asn$ on which $Q$ acts invertibly. In order to
describe these modules, the following notation is helpful.

\definition{Definition \secb.7}
We define the {\it trivial module} for $\asn$ to be the 1-dimensional 
$k$-vector space on which all the $u_i$ act as zero.

If $N$ is a left $\kq$-module, we can define an induced $\asn$-module
$\asn \otimes_{\kq} N$, where $\asn$ is regarded as a right module in which
$q$ acts as (right) multiplication by $Q$.
In particular, if $c \in k$ and $m \in \enn$, we define 
the $\asn$-module 
$\asml{c}{m}$ to be $$
\asml{c}{m} := \asm \otimes_{\kq} \frac{\kq}{\langle (q - c)^m \rangle}
.$$
\enddefinition

The fact that the trivial module is well-defined follows by a routine check
of the relations in Definition \secb.1.

The next result, which we prove in Section \secf, 
describes the finite dimensional indecomposable modules for
$\asn$, except those for which $Q$ acts nilpotently. The classification of
the latter modules is provably hopeless, as we will see in 
Theorem \secd.11 (iii).

\proclaim{Theorem \secb.8}
Let $k$ be an algebraically closed field, let $m \geq 1$ be a natural number,
and let $c$ be a nonzero element of $k$.
\item{\rm (i)}{The $\asn$-module $\asml{c}{m}$ has dimension $2^nm$. It 
is indecomposable in general, and irreducible if $m = 1$.}
\item{\rm (ii)}{If $M$ is a finite dimensional irreducible
$\asn$-module, then either $M$ is the trivial module, or $M$ is isomorphic
to $\asml{c}{1}$ for a unique value of $c$.}
\item{\rm (iii)}{If $M$ is a finite dimensional indecomposable
$\asn$-module on which $Q$ acts invertibly, then $M$ is isomorphic 
to $\asml{c}{m}$ for unique values of $c$ and $m$.}
\endproclaim

\head \secba. Heaps and fully commutative elements \endhead

In Section \secba, we will review some properties of heaps that will be
needed in the sequel.

Let $W$ be a Coxeter group with generating set $S$.
If $w \in W$ is a fixed element, then it may be
expressed as a product $w = s_1 s_2 \cdots s_r$ of elements of $S$, where the
$s_i$ are not necessarily distinct.
If $r$ is minimal subject to $w$ equalling such a product, then one calls 
$r$ the {\it length} of $w$, denoted by $\ell(w)$, and the
expression $\bw = s_1 s_2 \cdots s_r$
is called a {\it reduced expression} for $w$. The {\it support} of $w$,
$\supp(w)$, is
the set of elements of $S$ that appear in a reduced expression of $w$; this
subset is independent of the expression chosen. If the support of $w$ is the
whole of $S$, then we say that $w$ has {\it full support}.

Two expressions for a word $\bw \in S^*$ in the generators $S$
are called {\it commutation equivalent}
if it is possible to transform one to the other by iterated commutation of
adjacent generators; that is, relations of the form $s_i s_j = s_j s_i$, where
$m(i, j) = 2$. The associated equivalence classes of expressions are called
{\it commutation classes}. If $w$ has a single commutation class of reduced
expressions, we call $w$ {\it fully commutative}. 

If $\bw = s_1 s_2 \cdots s_r$ is any word in $S$, and $(i_1, i_2, \ldots, i_t)$
is a sequence satisfying $1 \leq i_1 < i_2 < \cdots < i_t \leq r$, then we
call the word $\bw'= s_{i_1} s_{i_2} \cdots s_{i_t}$ a {\it subexpression} of
$\bw$. A consecutive subexpression of $\bw$ (that is, one of the form
$s_i s_{i+1} s_{i+2} \cdots s_j$) is known as a {\it subword} of $\bw$.

Fully commutative elements can be characterized in terms of subwords as 
follows.

\proclaim{Theorem \secba.1 (Stembridge \cite{{\bf 16}})}
Let $W$ be a Coxeter group and let $w \in W$. Then $w$ is fully commutative
if and only if no reduced expression for $w$ has an alternating subword
$s_i s_j s_i \cdots$ of length $m(i, j)$.
\endproclaim

Minuscule elements, which were introduced by Stembridge in \cite{{\bf 17}}, 
can be defined by more 
restrictive criteria. These require a Dynkin diagram
for their definition, but in this paper, we will restrict our attention to the
Dynkin diagrams of type $\twd$ $(n \geq 2)$, shown in Figure \secba.1.

\topcaption{Figure \secba.1} Dynkin diagram of type $\twd$
\endcaption

\hskip 1.5cm
\hbox to 271pt{
\vbox to 24pt{\vfill
        \includegraphics{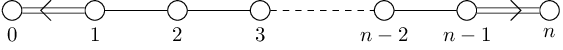}
}
\hfill}

A Dynkin diagram $\Gamma$ gives rise to a Coxeter group $W(\Gamma)$ 
(known as a Weyl group)
by associating a generator of $S$ to each vertex, and defining $$
m_{s, t} = \cases
1 & \quad \text{ if } s = t,\cr
2 & \quad \text{ if $s$ and $t$ are not adjacent,}\cr
3 & 
\quad \text{ if $s$ and $t$ are connected by an undecorated edge, and } \cr
4 & 
\quad \text{ if $s$ and $t$ are connected by a double edge with a single 
arrow.}\cr
\endcases
$$

By \cite{{\bf 17}, Proposition 2.3}, the following is a well-defined 
characterization of minuscule elements, which we take as our definition.

\definition{Definition \secba.2}
Let $\Gamma$ be a Dynkin diagram of type $\twd$ $(n \geq 2)$, and let
$W = W(\Gamma)$ be the associated Weyl group. Let $s_{i_1} s_{i_2} \cdots s_{i_r}$
be a word in the generators for $W$. For $1 \leq a < b \leq r$, we call a 
pair of generators $s_{i_a}$ and $s_{i_b}$ {\it consecutive occurrences of
$s_p$} if we have both $i_a = i_b = p$ and $i_c \ne p$ for $a < c < b$.

If $s_{i_1} s_{i_2} \cdots s_{i_r}$ is a reduced expression for $w$, then we 
call $w$ {\it minuscule} if for all labels $p$, whenever $s_{i_a}$ and 
$s_{i_b}$ are consecutive occurrences of $s_p$ with $1 \leq a < b \leq r$,
the sequence $I = (i_{a+1}, i_{a+2}, \ldots, i_{b-1})$ satisfies one of the 
following two conditions:
\item{(i)}{$I$ contains precisely two occurrences of labels adjacent to $p$,
and for each such label $q$, either $p$ and $q$ are connected by a single edge,
or $p$ and $q$ are connected by a double edge with a single arrow pointing
towards $q$; or}
\item{(ii)}{$I$ contains precisely one occurrence of a label $q$ 
adjacent to $p$, and $p$ and $q$ are connected by a double edge with an 
arrow pointing towards $p$.}

We denote the set of minuscule elements of $W$ by $W_m$.
\enddefinition

In general, there is no easy way to characterize minuscule elements in terms
of subword avoidance, but it is possible to do this
in the particular case of the Dynkin diagram $\twd$.
In the context of minuscule elements for $\twd$, we will call an expression
{\it forbidden} if either (a) it is of the form $s_i s_i$ or (b) it is of
the form $s_i s_j s_i$ where 
$i$ and $j$ are adjacent, and connected either 
by a single edge or by a double edge with an arrow pointing towards $j$.

\proclaim{Proposition \secba.3}
Let $W$ be a Coxeter group of type affine $C_n$ (identified with the Weyl
group of $\twd$), and let $\bw$ be a word in the generating set $S$.
\item{{\rm (i)}}{If $p$ and $q$ are adjacent labels and there exist consecutive
occurrences of $s_p$ in $\bw$
that are not separated by an occurrence of $s_q$,
then $\bw$ is commutation equivalent to an expression with a forbidden
subword.}
\item{{\rm (ii)}}{The expression $\bw$ is a reduced expression for a 
minuscule element for $\twd$ 
if and only if it is not commutation equivalent to a reduced expression with
a forbidden subword.}
\endproclaim

\demo{Proof}
We first prove (i). Suppose that $s_{i_a}$ and $s_{i_b}$ are consecutive 
occurrences
of $s_p$, where $p \not\in \{0, n\}$, but (without loss of generality) there
is no occurrence of $p-1$ in the sequence $I = (i_a + 1, \ldots, i_b - 1)$.
We assume that $a$ and $b$ are chosen so that $b-a$ is as small as possible.
We cannot have two (or more) occurrences of $p+1$ in the sequence $I$ because
there would be no occurrence of $s_p$ between them, contradicting minimality.
If there are no occurrences of $s_{p+1}$ in $I$ then $\bw$ is commutation 
equivalent to an expression containing the forbidden subword $s_p s_p$.
The other possibility is that we are in the situation of 
Definition \secba.2 (ii), with a
single occurrence of $p+1$ in the sequence $I$. This cannot happen either,
because $p$ is not an endpoint of $\Gamma$, which means that
the arrow between $p$ and $p+1$ in $\twd$, if there is one, 
points towards $p+1$. This in turn means that $\bw$ is commutation equivalent
to an expression containing the forbidden subword $s_p s_{p+1} s_p$.

We next prove (ii). No reduced expression can
have a subword of the form $s_p s_p$, and if one of the forbidden subwords 
of the form $s_i s_j s_i$ appears in some reduced
expression for $w$, then that reduced expression would violate the conditions
of Definition \secba.2. It follows that the given conditions are necessary
in order to have a reduced expression for a minuscule element.

Now suppose that $\bw \in S^*$ is not a reduced expression for a 
minuscule element; we will be done if 
can show that $\bw$ is commutation equivalent to an expression containing 
a forbidden subword.  Let $i_a$, $i_b$, $I$ and $p$ be chosen so as to 
violate one of the conditions of Definition \secba.2, and let $l$ be the 
number of labels in $I$ that are adjacent to $p$. 

If $l = 0$, then $\bw$ is commutation equivalent to an expression containing
the forbidden subword $s_p s_p$.
If $l = 1$, $\bw$ is commutation equivalent to a word with a subword 
$s_p s_q s_p$, and we must have violated condition
(ii) of Definition \secba.2. In this case, the hypotheses about the edge between
$p$ and $q$ show that $s_p s_q s_p$ is a forbidden subword.

If $l = 2$, we must have violated condition (i) of Definition \secba.2. 
Here, there are two elements of $I$ with labels adjacent to $p$. 
If the elements have distinct labels, $p$ is not an endpoint of $\Gamma$ 
the nature of the diagram $\twd$ means that
the hypotheses about the edges in Definition \secba.2 are automatically
satisfied. It must therefore be the case that these two elements of $I$ 
have the same label, and then we are done by part (i) of the current result. 

If $l \geq 3$, the fact that there are at most two labels adjacent to $p$ means
that there exist two elements of $I$ with label $q$, for the same neighbour
$q$ of $p$. These two occurrences of $q$ have no occurrence of $p$ between
them, and part (i) of the current result applies to complete the proof of (ii).
\qed\enddemo

The relevance of Proposition \secba.3 is that if $w$ is not a minuscule
element, and $s_{i_1} s_{i_2} \cdots s_{i_r}$ is a reduced expression for $w$, 
then it follows from Definition \secb.1 that the corresponding word
$u_{i_1} u_{i_2} \cdots u_{i_r}$ in $\asn$ is zero. In particular, the elements
$s_1 s_0 s_1$ and $s_{n-1} s_n s_{n-1}$ are not minuscule even though they are
fully commutative. On the other hand, the words $s_0 s_1 s_0$ and 
$s_n s_{n-1} s_n$ are minuscule (and therefore also fully commutative).

It will turn out (see Theorem \secd.11 (i)) 
that the minuscule elements of $W$ index a basis for $\asn$.

The Dynkin diagram of type affine $C_n$ can be obtained from that
of type $\twd$ by reversing the direction of the two arrows. The
two Dynkin diagrams give rise to isomorphic Weyl groups, but different sets of
minuscule elements. It will turn out that the representation theory of
$\asn$ is governed by the spin representation of the Lie algebra of type $B_n$,
which has dimension $2^n$, rather than that of type $C_n$. Because of this, a
case could be made that the elements in Definition \secba.2 should properly
be called ``cominuscule elements of type affine $C$''.
Similarly, the algebra $\asn$ could (more accurately) be called
``the nil Temperley--Lieb algebra of type twisted affine $D$'', but we choose
not to do this for reasons of brevity, since the underlying Weyl group
is of type affine $C$.

We now introduce heaps, which are certain labelled partially ordered sets 
that are closely related to fully commutative elements.
The treatment of heaps here follows that of \cite{{\bf 10}}.

\definition{Definition \secba.4}
A {\it heap} is a function $\e: E \rightarrow \Gamma$, where $E$
is a poset and $\Gamma$ is a graph, satisfying the following two
conditions.
\item{(i)}{The inverse images of each vertex $\e^{-1}(a)$ and each edge
$\e^{-1}(\{a, b\})$ are chains in $E$. Such chains are known as {\it vertex
chains} and {\it edge chains}, respectively.}
\item{(ii)}{The partial order $\leq$ on $E$ is the smallest partial order in
which the subsets in (i) above are chains.}
\enddefinition

A particularly important heap for our purposes arises from the following
poset.

\definition{Definition \secba.5}
Let $n \geq 2$, and let $\heapel$ be the set $$
\heapel = \{(a, b) \in \zed \times \zed: 0 \leq a \leq n \text{\ and\ }
a-b \text{\ is\ even}\}
.$$ 

We define the relation $\prec_\heapel$ on $\heapel$ by the condition 
$(a, b) \prec (c, d)$ if and only if $c = a \pm 1$ and $d = b + 1$. 
The relation $\leq$ on $\heapel$ is defined to be the reflexive, transitive
extension of $\prec$.
We may denote $\prec_\heapel$ by $\prec$ for short if the context is clear.
\enddefinition

Recall that a partially ordered set is called {\it locally finite} if all of
its intervals are finite. The partial order on a locally finite poset is always
the reflexive, transitive closure of its covering relations.
If $E$ is locally finite, we call $E$ {\it alternating} if no edge chain of
$E$ contains two successive vertices with the same label.

\proclaim{Lemma \secba.6}
With the notation above, $(\heapel, \leq)$ is a locally finite,
partially ordered set with covering relations given by $\prec$.
The function $\e : \heapel \ra \Gamma = \twd$ given by 
$\e((a, b)) = a$ makes $\heapel$ into an alternating heap.
\endproclaim

\demo{Proof}
The assertion about $E$ being alternating follows from the fact that the 
covering pairs in $E$ have adjacent labels. The other assertions are a
restatement of \cite{{\bf 10}, Proposition 6.4.14}.
\qed\enddemo

\remark{Remark \secba.7}
The partial order on $\heapel$ may be defined directly as follows:
$(a, b) \leq (c, d)$ if and only if both $b \leq d$ and $|c-a| \leq |d-b|$
(see \cite{{\bf 10}, Definition 6.1.1}). Note that, because $a \equiv b \mod 2$
and $c \equiv d \mod 2$, it follows that $|c-a| \equiv |d-b| \mod 2$.
\endremark

A heap can be depicted in terms of its labelled Hasse diagram; that is, 
a Hasse diagram in which each element $\alpha \in E$ is labelled by
$\e(\alpha) \in \Gamma$. Figure \secba.2 shows the labelled Hasse diagram 
for the heap $\heapel$ in the special case $n = 6$.

\topcaption{Figure \secba.2} The heap $E(6)$
\endcaption

\hskip 1.5cm
\hbox to 233pt{
\vbox to 171pt{\vfill
        \includegraphics{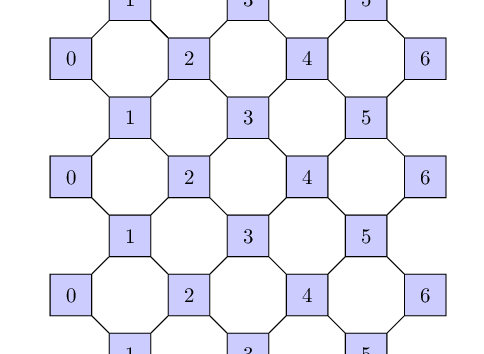}
}
\hfill}

Heaps can be made into the objects of a category, as follows.

\definition{Definition \secba.8}
There is a category $\hcat$ whose objects are heaps, in which a 
morphism $f$ from
a heap $(E_1, \leq_1)$ over a graph $\Gamma_1$ to a heap $(E_2, \leq_2)$
over a graph $\Gamma_2$ consists of a pair $(f_E, f_\Gamma)$ in which 
\item{(i)}{$f_E$ is a morphism of partially ordered sets 
(i.e., $x \leq y \Rightarrow f_E(x) \leq f_E(y)$);}
\item{(ii)}{$f_\Gamma$ is a morphism
of graphs (i.e., if $a$ and $b$ are adjacent, then the vertices $f_\Gamma(a)$
and $f_\Gamma(b)$ are adjacent or equal) and}
\item{(iii)}{the following diagram commutes:

\centerline{
\hbox to 66pt{
\vbox to 56pt{\vfill
        \includegraphics{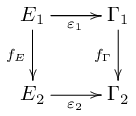}
}
\hfill}
}}

The category $\hcat$ has a subcategory, $\hcatg$,
whose objects are the heaps over $\Gamma$.  A 
morphism $f$ of $\hcatg$ is
a morphism $f \in \Hom_\hcat(A, B)$, where $A$ and $B$ are objects of
$\hcatg$, and where $f_\Gamma$ is the identity map.
\enddefinition

For our purposes, $\Gamma$ will be the Dynkin diagram $\twd$, and
we will work with the category $\hcatg$. The isomorphisms in this category
are label-preserving isomorphisms of partially ordered sets.

\definition{Definition \secba.9}
Let $E_1$ and $E_2$ be heaps over $\Gamma$.  We say that $E_2$ is a 
{\it subheap}
of $E_1$ if there exists a morphism $$f \in \Hom_\hcatg(E_2, E_1)$$ that is
injective on vertices (i.e., $f_E$ is injective).  We will often 
identify $E_2$ with the subset of
$E_1$ given by the image of $f$ on vertices.  

We say that the subheap $E_2$ is an {\it ideal} of $E_1$ if 
whenever $y \in E_2$ and $x \in E_1$ satisfy $x \leq y$, then we have 
$x \in E_2$. 
Dually, we say that the subheap $E_2$ is a {\it filter} of $E_1$ if 
whenever $y \in E_2$ and $x \in E_1$ satisfy $x \geq y$, then we have 
$x \in E_2$. 
We say that the subheap $E_2$ is
{\it convex} if it corresponds to a convex subset of $E_1$, meaning that 
whenever $\alpha < \beta < \gamma$ are elements of $E_1$ with $\alpha, \gamma
\in E_2$, then $\beta \in E_2$. It is immediate from the definitions that
ideals and filters are convex.

We say that $E_2$ is a {\it chain} if it is a chain in $E_1$ considered 
as a partially ordered set.

We say that $E_2$ is a {\it closed $p$-interval} if it is an interval 
of the form $[x, y]$, where $x, y \in E$ satisfy $\e(x) = \e(y) = p$ and 
$(x, y) \cap \e^{-1}(p) = \emptyset$. If $[x, y]$ is a closed $p$-interval,
we call $(x, y)$ an {\it open $p$-interval}.
\enddefinition

The shaded area in Figure \secba.3 shows a subheap $F$ of the heap
$E(6)$. In this case, $F$ has nine elements and is convex, but is not a chain.

\topcaption{Figure \secba.3} A finite subheap of $E(6)$
\endcaption

\hskip 1.5cm
\hbox to 233pt{
\vbox to 171pt{\vfill
        \includegraphics{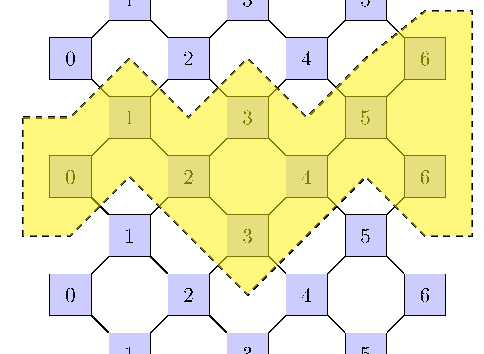}
}
\hfill}

\remark{Remark \secba.10}
If $F$ and $E$ are heaps over $\Gamma$, and $F$ is a subheap of $E$, the
partial order on $F$ may not agree with the restriction to $F$ of the partial
order on $E$. Consider, for example, the subheap $F$ of $E(2)$ consisting of
the two elements $\alpha = (0, 0)$ and $\beta = (2, 2)$. We then have
$\alpha < \beta$ as elements of $E(2)$, and yet the partial order on $F$ is
the trivial one.

However, if $F$ is a convex subheap of $E$, the partial orders
will be compatible in this way (see \cite{{\bf 10}, Exercise 2.1.7}). The 
incompatibility in the example above arises from the fact that 
$\{(0, 0), (2, 2)\}$ is not a convex subset of $E(2)$: we have 
$(0, 0) < (1, 1) < (2, 2)$, but $(1, 1)$ does not lie in $F$.
\endremark

A key property of finite heaps is that they correspond to elements of the
commutation monoid, which is defined as follows.

\definition{Definition \secba.11}
Let $\Gamma$ be a Dynkin diagram with vertex set $S$. The
{\it commutation monoid}, $\cog$, of 
$\Gamma$ is the quotient of the free monoid $S^*$ by the congruence 
$\equiv$ generated by the relations $st \equiv ts$ whenever $s$ and $t$ 
are nonadjacent in $\Gamma$.  The equivalence class of the word $x \in S^*$
is denoted by $[x]$.  The multiplication in $\cog$ is given by $$
[x][y] := [xy]
,$$ and is well-defined.
The elements of $\cog$ are called {\it traces}. 
\enddefinition

A trace is an equivalence class of words in the alphabet $S$, up to 
commutation. It is immediate from the definitions that the Weyl group $W$ 
is a quotient of the commutation monoid, and that the algebra $\asn$ is a 
quotient of the monoid algebra of this monoid. It follows that the reduced
expressions for an element $w \in W$ correspond to a single trace if and only
if $w$ is fully commutative, and that if $x = s_{i_1} s_{i_2} 
\cdots s_{i_r}$, then the trace $[x]$ corresponds to a well-defined element 
$u_x = u_{[x]} \in \asn$ given by $u_{[x]} = u_{i_1} u_{i_2} \cdots u_{i_r}$.

Suppose $(F, \leq_F)$ is a finite heap over $\Gamma$. The partial order 
$\leq_F$ can be extended to a total order, $\leq'$ on $F$ (by Szpilrajn's
Theorem). Let $\a_1, \a_2, \ldots, \a_k$ be the elements of $F$, listed in
the order $\a_1 >' \a_2 >' \cdots >' \a_k$. Let $f_{>'}(F)$ be the
element of $\cog$ given by $\e(\a_1)\e(\a_2)\cdots \e(\a_k)$.

\proclaim{Proposition \secba.12}
The trace $f(F) = f_{>'}(F)$ is independent of the choice of total order 
$\leq'$. The resulting function $f$ gives a bijection between the isomorphism
classes of finite heaps over $\Gamma$ and the elements of $\cog$.
\endproclaim

\demo{Proof}
This result is essentially due to Viennot \cite{{\bf 18}, Proposition 3.4}. It 
is stated in the above notation (but with the opposite total order) 
in \cite{{\bf 10}, Theorem 2.1.20}.
\qed\enddemo

In the notation of Proposition \secba.12, we will say that the heap $F$
{\it represents} the trace $f(F)$. Because the reduced expressions of a fully
commutative element correspond to the same element of $\cog$, it follows
that each fully commutative element of $W$ has a unique heap, up to
isomorphism in $\hcatg.$ We will call this heap $F$ {\it the heap of $w$} for 
short, and we will also say that $F$ represents $w$. 
For example, Figure \secba.3 shows the heap of the minuscule element
$s_6 s_1 s_3 s_5 s_0 s_2 s_4 s_6 s_3$. 

As noted in the previous proof, we use the opposite ordering from the 
conventional one when associating heaps with words. The reason for this is
that we will be using left modules rather than right modules, and raising 
operators rather than lowering operators.

\proclaim{Lemma \secba.13}
Suppose that $w \in W$ is a minuscule element of the form $w = uv$, where
$\ell(w) = \ell(u) + \ell(v)$. Then the heap $F$ of $w$ can be expressed 
as a disjoint union of the form $F = F_1 \cup F_2$, where $F_1$ is the heap of
$u$, $F_2$ is the heap of $v$, $F_1$ is a filter of $F$, and $F_2$ is an ideal
of $F$.
\endproclaim

\demo{Proof}
Let $s_{i_1} s_{i_2} \cdots s_{i_l}$ be a reduced expression for $u$, and
let $s_{j_1} s_{j_2} \cdots s_{i_m}$ be a reduced expression for $v$. By 
hypothesis, $s_{i_1} \cdots s_{i_l} s_{j_1} \cdots s_{j_m}$ is a reduced
expression for $w$, meaning that there is a reduced expression for $w$ in
which none of the generators arising from $v$ appear to the left of any
of the generators arising from $u$. The result follows from the correspondence
of Proposition \secba.12.
\qed\enddemo

\head \secc. Convex subheaps \endhead

The main purpose of Section \secc\ is to prove Theorem \secc.13, which shows
that any minuscule element can be represented by a convex subheap of 
$\heapel$. As we mention in the conclusion, this is closely related to
a certain property of Coxeter elements. A very useful ingredient in the 
proof is Theorem \secc.1, which shows how to characterize a convex subheap
with full support using properties of the edge chains. This result is general
and appears to be new.

\proclaim{Theorem \secc.1}
Let $\Gamma$ be a graph with no isolated points, let $\e : E \ra \Gamma$ be
a heap, and let $F$ be a subheap of $E$ with full support. If, for every 
edge $\{s, t\}$ of $\Gamma$, the edge chain $F \cap \e^{-1}(\{s, t\})$
is convex as a subposet of the edge chain $E \cap \e^{-1}(\{s, t\})$, then 
$F$ is convex as a subheap of $E$.
\endproclaim

\demo{Proof}
Since $\Gamma$ has no isolated points, every vertex of $\Gamma$ is part of
an edge, and every vertex chain of $E$ is a subset of an edge chain of $E$.
(In particular, this means that any element of $E$ is an element of some
edge chain.)
In general, the partial order on a heap is the reflexive, transitive extension 
of the vertex and edge chain relations, but under our hypotheses, the 
partial order on $E$ can be recovered from the edge chains alone.

Let $I$ be the subset of $E$ consisting of all elements $x \in E$ for
which there exists $y \in F$ (depending on $x$) 
satisfying $x \leq y$ and $\e(x) = \e(y)$. It follows immediately that
we have $F \subseteq I$, and we claim that $I$ is an ideal of
$E$. By the above paragraph, it is enough to show that whenever $x$ and $z$ are
in the same edge chain and satisfy $x \in I$ and $z \leq x$, then we have
$z \in I$.

Assume that $x$ and $z$ have the above properties, and that they lie in
the same edge chain, $E_{s, t} = E \cap \e^{-1}(\{s, t\})$, of $E$. Without 
loss of generality, we may assume that $s = \e(x)$. By definition of $I$,
there exists $y \in F \cap E_{s, t}$ such that $x \leq y$.
If there exists an element $v \in F \cap E_{s, t}$ with $v \leq z$, then we
can apply the convexity assumption on $F \cap E_{s, t}$ to the triple 
$v \leq z \leq y$ to deduce that $z \in F$, which proves that $z \in I$, as 
required.

The other possibility is that we have $z < v$ for all elements $v$ in
the chain $F \cap E_{s, t}$. Since $F$ has full support, there is an element
$y' \in F$ for which $\e(y') = \e(z)$, which implies that 
$y' \in F \cap E_{s, t}$. By assumption, we have $z < y'$, and the definition
of $I$ (applied to $z$ and $y'$) now shows that $z \in I$.

The argument above shows that the edge chain $I \cap E_{s, t}$ consists
of the disjoint union of (a) the set $F \cap E_{s, t}$ with (b) the set of 
all elements of $E_{s, t}$ that are strictly less than all elements of
$F \cap E_{s, t}$.

Dually, we define $I'$ be the subset of $E$ consisting of all elements 
$x \in E$ for which there exists $y \in F$ (depending on $x$) 
satisfying $x \geq y$ and $\e(x) = \e(y)$. Arguing as above, we find that
$I'$ is a filter of $E$ with the property that each edge chain
$I' \cap E_{s, t}$ consists
of the disjoint union of (a) the set $F \cap E_{s, t}$ with (b) the set of 
all elements of $E_{s, t}$ that are strictly greater than all elements of
$F \cap E_{s, t}$.

Since ideals and filters are always convex, and intersections of convex
subsets are convex, it follows that $F' = I \cap I'$ is convex. The
previous two paragraphs show that for any edge chain $E_{s, t}$, we
have $F' \cap E_{s, t} = F \cap E_{s, t}$. It follows that $F' = F$, and thus
that $F$ is a convex subheap of $E$.
\qed\enddemo

\remark{Remark \secc.2}
It follows from the definitions that any convex subheap $F$ of a heap $E$
must satisfy the condition that each of its edge chains 
$F \cap \e^{-1}(\{s, t\})$ is convex as a subset of $E \cap \e^{-1}(\{s, t\})$.
However, the hypothesis on the support of $F$ cannot be removed, as
the counterexample in Remark \secba.10 shows.
\endremark

Unless otherwise stated, the results in the rest of this section
only apply in the context of the Dynkin diagram $\Gamma = \twd$. In
particular, the Weyl group $W$ refers to $W(\twd)$, and each heap $F$ is a heap
over $\twd$.

\proclaim{Lemma \secc.3}
Let $F$ be an alternating heap over $\twd$, and let 
$x \prec y$ be a covering relation in $F$. Then $\e(x)$ and $\e(y)$ are
adjacent labels.
\endproclaim

\demo{Proof}
Since $p = \e(x)$ is not an isolated point in $\Gamma$, there exists a label
$q$ that is adjacent to it.
The definition of the partial order on $F$ means that $x$ and $y$
are either part of the same vertex chain, or part of the same edge chain.
The former possibility cannot occur, because $x$ and $y$ would be adjacent
elements with the same label of the edge chain $F \cap \e^{-1}(\{p, q\})$, 
contradicting the alternating condition.
\qed\enddemo

The next result shows that the edge chains of the heap of a minuscule element
consist entirely of covering relations. (This is not true for minuscule
elements in all types, and there are counterexamples if the Dynkin diagram 
contains a circuit.)

\proclaim{Lemma \secc.4}
If $F$ is an alternating heap in type $\twd$, then 
each edge chain in $F$ consists entirely of covering relations (in $F$).
\endproclaim

\demo{Proof}
Let $x < y$ be an adjacent pair of elements in an edge chain $F \cap
\e^{-1}(\{s, t\})$. We may assume without loss of 
generality that $\e(x) = s$ and $\e(y) = t$. It remains to show that 
$x \prec y$ is a covering relation, so suppose that this is not the case.

Since we have $x < y$, there exists a chain in $F$ of covering relations of
the form $$
x = \alpha_1 \prec \alpha_2 \prec \cdots \prec \alpha_k = y
.$$ By Lemma \secc.3, the labels $\e(\alpha_i)$ and $\e(\alpha_{i+1})$ are
adjacent for each $1 \leq i < k$. The Dynkin diagram $\Gamma$ has no circuits,
so every path from $\e(x)$ to $\e(y)$ must traverse the edge $\{s, t\}$.
It follows that for some $i$, we must have $\{\e(\alpha_i), \e(\alpha_j)\}
= \{s, t\}$, but the assumptions on $x$ and $y$ show that they are the only
elements in the sequence of $\alpha_i$ with labels $s$ or $t$. We conclude
that $k = 2$ and that $x \prec y$ is a covering relation, as required.
\qed\enddemo

Heaps of minuscule elements may be characterized as follows in terms of their 
$p$-intervals: each open $p$-interval contains 
exactly one element with each label adjacent to $p$, and no other elements.
More precisely, we have

\proclaim{Proposition \secc.5}
Let $W$ be the Weyl group associated to a Dynkin diagram 
$\Gamma = \twd$ for $n \geq 2$, and let 
$F$ be a heap representing an expression $\bw \in S^*$.
Then $\bw$ is a reduced expression for a minuscule element 
if and only if every open $p$-interval $(x, y)$ of $F$ 
satisfies one of the following conditions:
\item{\rm (i)}{$p = 0$ and $(x, y) = \{z\}$, where $\e(z) = 1$;}
\item{\rm (ii)}{$0 < p < n$ and $(x, y) = \{z_1, z_2\}$, where 
$\e(z_1) = p-1$ and $\e(z_2) = p+1$;}
\item{\rm (iii)}{$p = n$ and $(x, y) = \{z\}$, where $\e(z) = n-1$.}
\endproclaim

\demo{Proof}
If $p = 0$ or $p = n$, then Definition \secba.2 (ii) shows that $(x, y)$ 
contains
a single element $z$ for which $q = \e(z)$ is adjacent to $p$; in particular,
we have $q = 1$ if $p = 0$, and $q = n-1$ if $p = 1$. It follows
that $x$, $z$ and $y$ are consecutive elements of the edge chain 
$\e^{-1}(\{p, q\})$, and Lemma \secc.4 shows that $x \prec z$ and $z \prec y$
are covering relations. Lemma \secc.3 then shows that $(x, y)$ contains
no other elements, because no element other than $z$ can cover $x$ or be
covered by $y$.

If $0 < p < n$ then Definition \secba.2 (i) shows that $(x, y)$ contains
precisely two elements $z_1$ and $z_2$ with labels adjacent to $p$, and 
Proposition \secba.3 (i) shows that, without loss of generality, we have
$\e(z_1) = p-1$ and $\e(z_2) = p+1$. Considering the edge chains 
$\e^{-1}(\{p-1, p\})$ and $\e^{-1}(\{p, p+1\})$, and arguing as in the previous
paragraph, we find from Lemma \secc.4 that $x \prec z_i$ and $z_i \prec y$ 
are covering relations for $i \in \{1, 2\}$. Again, Lemma \secc.3 then shows
that $(x, y)$ contains no other elements.
\qed\enddemo

Recall that if $E$ is a poset, a function $\rho : E \ra \zed$ is said to be a
{\it rank function} for $E$ if whenever $a, b \in E$ are such that
$a \prec b$ is a covering relation, we have $\rho(b) = \rho(a) + 1$.  It 
follows that if $E$ is locally finite and $x, y \in E$ satisfy $x < y$,
then $\rho(x) < \rho(y)$. If a
rank function for $(E, \leq)$ exists, we say $(E, \leq)$ is {\it ranked}.

\proclaim{Lemma \secc.6}
The heap of a minuscule element $w \in W$ is alternating and ranked.
\endproclaim

\demo{Proof}
Denote the heap of $w$ by $F$, and let $s$ and $t$ be adjacent labels.
By Proposition \secba.3, if $x, y \in F$ satisfy $\e(x) = \e(y) = s$, then
there must exist $z \in F$ with $x < z < y$ and $\e(z) = t$. In other words,
any two vertices in $F$ with label $s$ must be
separated by a vertex with label $t$, and vice versa; this proves that $F$
is alternating.

By \cite{{\bf 9}, Theorem 2.1.1}, a finite heap over a graph with no circuits 
is ranked if and only if all of its closed $p$-intervals are ranked. The latter
condition follows from the characterization of $p$-intervals given in
Proposition \secc.5.
\qed\enddemo

A rank function can be used to identify the heap of a minuscule element with 
a subheap of $\heapel$, using the following embedding.

\proclaim{Lemma \secc.7}
Let $F$ be the heap of a minuscule element $w \in W$, and let 
$\rho : F \ra \zed$ be a rank function for $F$.
Define a map $\iota_\rho : F \ra \heapel$
by $\iota_\rho(x) = (\e(x), \rho(x))$.
\item{\rm (i)}{The map $\iota_\rho$ sends covering relations
in $F$ to covering relations in $\heapel$.}
\item{\rm (ii)}{The map $\iota_\rho$ is a morphism in $\hcatg$
that is injective on vertices, so that the image of $\iota_\rho$ is
a subheap of $\heapel$.}
\endproclaim

\demo{Proof}
It is immediate from the definitions that $\iota_\rho$ sends elements of $F$
to elements of $\heapel$ with the same label. Furthermore, if 
$x, y \in F \cap \e^{-1}(p)$ satisfy $x < y$,
then we have $\rho(x) < \rho(y)$, and in turn
$\iota_\rho(x) \ne \iota_\rho(y)$. It follows that $\iota_\rho$ is injective
on vertex chains, which completes the proof that $\iota_\rho$ is injective
on vertices.

Let $x, y \in F$ be such that $x \prec y$ is a covering
relation. Lemmas \secc.3 and \secc.6 show that $\e(x) = \e(y) \pm 1$, and the 
definition of $\iota_\rho$ shows that $\iota_\rho(x)$ and $\iota_\rho(y)$
are of the form $(a, b)$ and $(a \pm 1, b+1)$, respectively. The latter is
a covering relation in $\heapel$ by definition, which proves (i).

Since $F$ and $\heapel$ are
both locally finite, each partial order is determined by its covering
relations, and it follows that $\iota_\rho$ is a morphism of posets.
Because $\iota_\rho$ also respects labels, it follows that $\iota_\rho(F)$ 
is a morphism in $\hcatg$. The definition of subheap then assures that
$\iota_\rho(F)$ is a subheap of $\heapel$, which completes the proof of
(ii).
\qed\enddemo

\proclaim{Lemma \secc.8}
Let $w \in W$ be a minuscule element with full support, let $F$ be the heap of
$w$ and let $\rho : F \ra \zed$ be a rank function for $F$. Then 
$\iota_\rho(F)$ is a convex subheap of $\heapel$, and therefore the partial
order on $\iota_\rho(F)$ induced by $\heapel$ agrees with the one inherited
from $F$.
\endproclaim

\demo{Proof}
The assertion about inherited partial orders is a general property of
convex subheaps (see \cite{{\bf 10}, Exercise 2.1.7}).

Since $F$ has full support, to show $F$ is ranked it suffices 
by Theorem \secc.1 to show that
each edge chain of $\iota_\rho(F)$ is a convex subset of the corresponding 
edge chain for $\heapel$. 

By Lemmas \secc.4 and \secc.6, each edge chain
$F_{s, t} = F \cap \e^{-1}(\{s, t\})$ consists entirely of covering 
relations in $F$. Lemma \secc.6 also shows that $F$ is ranked, which means
that the elements of the chain $F_{s, t}$ have consecutive ranks.
Using the definition of $\iota_\rho$, we find that the corresponding edge chain
$\iota_\rho(F) \cap \e^{-1}(\{s, t\})$ consists of the set $$
\{(x, c) \in \heapel : \e(x) \in \{s, t\} \text{\ and\ } N_1 \leq c \leq N_2\}
,$$ for suitable integers $N_1 \leq N_2$. This shows that 
$\iota_\rho(F) \cap \e^{-1}(\{s, t\})$ is a convex as a subset of
$\heapel \cap \e^{-1}(\{s, t\})$, and this completes the proof.
\qed\enddemo

\definition{Definition \secc.9}
Let $F$ be a ranked heap with rank function $\rho$, over an arbitrary
graph $\Gamma$. 
If $F(p) = F \cap \e^{-1}(p)$ is nonempty and bounded above, let
$\rho_H(F, p)$ be the rank of the unique maximal element of $F(p)$.
If $F(p) = F \cap \e^{-1}(p)$ is nonempty and bounded below, let
$\rho_L(F, p)$ be the rank of the unique minimal element of $F(p)$.
\enddefinition

\proclaim{Lemma \secc.10}
Let $F$ be a ranked, alternating heap over $\twd$, with rank function $\rho$,
and let $p$ and $q$ be adjacent vertices of $\twd$.
Maintain the notation of Definition \secc.9.
\item{\rm (i)}{If $F(p)$ is nonempty and $F(q)$ is empty, then $F(p)$
consists of a single element.}
\item{\rm (ii)}{If $F(p)$ and $F(q)$ are both nonempty and bounded below,
then we have $\rho_L(F, p) = \rho_L(F, q) \pm 1$.}
\item{\rm (iii)}{If $F(p)$ and $F(q)$ are both nonempty and bounded above,
then we have $\rho_H(F, p) = \rho_H(F, q) \pm 1$.}

In particular, all these conditions hold if $F$ is the heap of a minuscule
element $w$.
\endproclaim

\demo{Proof}
If the edge chain $F \cap \e^{-1}(\{p, q\})$ of $F$ contains no elements
labelled $q$, the alternating condition means that the chain must
consist of a single element labelled $p$, proving (i).

Now suppose that the vertex chains $F(p)$ and $F(q)$ are both nonempty.
Suppose that $\alpha$ and $\beta$ be consecutive elements of 
the edge chain $F \cap \e^{-1}(\{p, q\})$ with $\alpha < \beta$.
By Lemma \secc.4, $\alpha \prec \beta$ is a covering relation, which
implies that $\rho(\beta) = \rho(\alpha) + 1$. 

If $F(p)$ and $F(q)$ are both bounded above, then so is the edge chain
$F \cap \e^{-1}(\{p, q\})$. If we take $\alpha$ and $\beta$ to be the highest
two elements of this edge chain, then (iii) follows.
A parallel argument establishes (ii).

The heap of a minuscule element is alternating and ranked by 
Lemma \secc.6, which proves the final assertion.
\qed\enddemo

\proclaim{Lemma \secc.11}
Let $w \in W$ be a minuscule element without full support, and let $u$
be the product of the generators in $S \backslash \supp(w)$, once each, but
in any order. Then the element $wu$ is a minuscule element with full support,
and we have $\ell(wu) = \ell(w) + \ell(u)$.
\endproclaim

\demo{Proof}
Let $\bw = s_{i_1} \cdots s_{i_k}$ be a reduced expression for $w$, and let
$\bu = s_{j_1} \cdots s_{j_l}$ be a reduced expression for $u$. We need to show 
that $\bx = \bw \bu$
is a reduced expression for a minuscule element; the assertion about full
support will then follow by the construction of $u$.

Suppose for a contradiction that this is not the case, which implies that
we must have two occurrences of a generator $s_p$ in $\bx$ that violate
the condition on $p$-intervals given in Proposition \secc.5. However, the
generator $s_p$ cannot occur in $\bu$, because 
the construction of $u$ ensures that each generator in $\supp(u)$ occurs 
only once in $\bx$. It follows that both these occurrences of $s_p$ occur
in $\bw$. 
Since the expression $\bw$ satisfies Proposition \secc.5,
we conclude that this situation cannot occur, and this completes the proof.
\qed\enddemo

\proclaim{Proposition \secc.12}
Each minuscule element $w \in W$ can be represented by a convex subheap of
$\heapel$.
\endproclaim

\demo{Proof}
If $w$ has full support, then the result follows from Lemma \secc.8, so 
we will assume that this is not the case. 

By Lemma \secc.11, there exists a minuscule
element $wu$ with full support, for which $\ell(wu) = \ell(w) + \ell(u)$.
Let $\bw$ and $\bu$ be reduced expressions for $w$ and $u$ respectively, and
let $F$ be the heap of $wu$. By Lemma \secc.8, there exists a rank function
$\rho : F \ra \zed$ such that $\iota_\rho(F)$ is a convex subheap of $\heapel$.

Define $F_w = F \cap \e^{-1}(\supp(w))$, and $F_u = F \cap \e^{-1}(\supp(u))$.
By Lemma \secba.13, $F_w$ is a filter of $F$ and $F_u$ is an ideal of $F$.
Since ideals and filters of convex subsets are themselves convex, we see
that the subheap $F_w$ satisfies the required conditions.
\qed\enddemo

We will make extensive use of the following result in later sections.

\proclaim{Theorem \secc.13}
Let $W = W(\Gamma)$ be the Weyl group corresponding to the Dynkin diagram of 
type $\Gamma = \twd$ $(n \geq 2)$, and let $w \in W$ be a minuscule 
element. Then any finite convex subheap of $\heapel$ is the heap of a minuscule
element of $W$, and any minuscule element can be represented by a finite
convex subheap of $\heapel$.
\endproclaim

\demo{Proof}
Proposition \secc.12 proves the second assertion.

To prove the first assertion, we note that the structure of the Hasse
diagram of $\heapel$ shows that the $p$-intervals of $\heapel$ satisfy the
conditions of Proposition \secc.5, and these conditions are inherited by
all convex subsets of $\heapel$. Any finite convex subheap of $\heapel$
corresponds to an expression in $S^*$ by Proposition \secba.12, and the
first assertion follows.
\qed\enddemo

The hypotheses of Theorem \secc.13 are more delicate than it might appear
at first, and the result fails for fully commutative elements in
general. For example, $w = s_1 s_0 s_1$ is a fully commutative element, and it
can be represented by a subheap of $\heapel$, but because any such subheap
contains an element labelled $2$ between any two elements labelled $1$, it
is not possible to represent $w$ by a convex subheap. 

Another subtlety is
that is possible for a minuscule element to be represented by a non-convex
subheap of $\heapel$. For example, the non-convex subheap $F$ in 
Remark \secba.10 
represents the minuscule element $x = s_0 s_2$. The element 
$x$ is also represented by the convex subheap of $\heapel$ given by 
$F' = \{(0, 0), (2, 0)\}$.
The heaps $F$ and $F'$ are isomorphic in $\hcatg$, even though one is a 
convex subheap of $\heapel$ and the other is not.

\head \secd. Weights \endhead

The main purpose of Section \secd\ is to prove Theorem \secb.4. In order
to do this, we will use a description of $\asm$ in terms of heaps.

We define an ideal $I$ of $\heapel$ to be a {\it proper ideal} if each 
vertex chain of $I$ satisfies $$
\emptyset \ne I \cap \e^{-1}(p) \ne \heapel \cap \e^{-1}(p)
.$$ Since $\heapel$ is a full heap (as defined in \cite{{\bf 10}, \S2.2}) and 
$\Gamma$ is finite and connected, it follows 
from \cite{{\bf 10}, Lemma 3.2.4 (v)}
that every ideal of $\heapel$ is proper with the exception of
$\heapel$ itself and the empty ideal. We denote the set of all proper ideals
of $\heapel$ by $\propid$.

\proclaim{Proposition \secd.1}
Let $k$ be a field, and let $V = V_\propid$ be the $k$-vector space with basis
$\{b_I : I \in \propid\}$. Then $V$
becomes a $\asn$-module in which each generator $u_i$ acts as a raising
operator, as follows: $$
u_i(b_I) = \cases
b_{I'} & \text{\ if\ } I'\backslash I = \{x\} \text{\ with\ } \e(x) = i;\cr
0 & \text{\ otherwise.}\cr
\endcases
$$
\endproclaim

\demo{Note}
If an element $x \in \heapel$ satisfies the conditions in the statement of
Proposition
\secd.1, it must be the unique minimal element of label $i$ in the
filter $\heapel\backslash I$. It follows from this observation
that the ideal $I'$ in the statement is unique if it exists.
\enddemo

\demo{Proof}
The raising operators $X_i$, corresponding to the action of the $u_i$,
are defined in \cite{{\bf 10}, Definition 3.1.4}. They
satisfy the defining relations of Definition \secb.1 by equations
(4.13)--(4.15) of \cite{{\bf 10}, Lemma 4.1.4}.
\qed\enddemo

A {\it Coxeter element} is an element $w \in W$ that can be written as the
product of the elements of the generating set $S$, once each, in some
order. Since the heap of a Coxeter element has no $p$-intervals, it follows
vacuously from Proposition \secc.5 that a Coxeter element is minuscule.
A Coxeter element $w$, by construction, has full support, so it can be
embedded into $\heapel$ as a convex subheap $F$ by using a rank function 
$\rho$, as in Lemma \secc.8. Lemma \secc.10 now shows that if $x$ and $y$ 
are the unique elements of $F$ with labels $p$ and $p+1$ respectively, then
we have $\rho(y) = \rho(x) \pm 1$. This shows that the two cases in the 
following definition are exclusive and exhaustive.

\definition{Definition \secd.2}
Let $F$ be a convex subheap of $\heapel$ representing a Coxeter element $w$.
We define the {\it contour} of $F$ to be the word $p_1 p_2 \cdots p_n$ in
the two-letter alphabet $\{{+}, {-}\}$, where we define $$
p_i = \cases
{+} & \quad \text{ if } \rho_L(F, i) = \rho_L(F, i-1) + 1, \text{ and }\cr
{-} & \quad \text{ if } \rho_L(F, i) = \rho_L(F, i-1) - 1.\cr
\endcases
$$ We will call the set of all $2^n$ contours {\it weights}, and we denote
the set of weights by $\Lambda$. 
\enddefinition

\remark{Remark \secd.3}
Because $F$ has a unique element of each label, we could just as well have used
$\rho_H$ instead of $\rho_L$ in Definition \secd.2. It is also possible to
define the contour solely in terms of a(ny) reduced expression $\bw$ for $w$, 
without reference to the heap: we have $p_i = {+}$ (respectively, $p_i = {-}$)
if $s_{i-1}$ occurs to the right (respectively, left) of $s_i$ in $\bw$.
Since Coxeter elements are minuscule and, therefore, fully commutative, it
then follows from Proposition \secba.12
that there is a bijection between contours and Coxeter elements. 
It also follows that there are $2^n$ Coxeter elements in the Weyl group of type
$\twd$.
\endremark

\definition{Definition \secd.4}
Let $F$ be a convex subheap of $E$ with full support. 
If each vertex chain of $F$ is nonempty and bounded above, then we define
$F_H$ to be the subheap of $F$ consisting of those elements of $F$ that 
are maximal in their vertex chains. 
If each vertex chain of $F$ is nonempty and bounded below, then we define
$F_L$ to be the subheap of $F$ consisting of those elements of $F$ that 
are minimal in their vertex chains. 
\enddefinition

\proclaim{Lemma \secd.5}
Let $F$ be a convex subheap of $\heapel$ with full support. 
\item{\rm (i)}{If $F$ has a subheap $F_L$ as in Definition \secd.4, then 
$F_L$ is an ideal of $F$, and $F_L$ is the heap of a Coxeter element of $W$.}
\item{\rm (ii)}{If $F$ has a subheap $F_H$ as in Definition \secd.4, then 
$F_H$ is a filter of $F$, and $F_H$ is the heap of a Coxeter element of $W$.}
\endproclaim

\demo{Proof}
Since $\heapel$ is locally finite and $F$ is a convex subheap,
the covering relations and partial order of $F$ are the restrictions of 
those of $\heapel$. Since $\heapel$ is ranked, it follows that any rank function
for $\heapel$ restricts to a rank function for $F$. We now fix a rank function,
$\rho$, for $F$.

Since $F$ is a convex subheap, each edge chain of $F$ is a convex subposet of
an edge chain of $\heapel$. Because $\heapel$ is alternating by Lemma \secba.6,
it follows that $F$ is alternating.

Let $x \in F$ be an element of $F_L$, and suppose that $\e(x) = p$. By
hypothesis, $x$ is the minimal element of the vertex chain $F \cap \e^{-1}(p)$.
We will be done if we can show that whenever $x$ covers an element $y \in F$,
then we have $y \in F_L$. By Lemma \secc.3, we have $\e(y) = q = p \pm 1$, and
the definition of rank function shows that $\rho(y) = \rho(x) - 1$, which
implies that $\rho_L(F, q) \leq \rho_L(F, p) - 1$.
Lemma \secc.10 (ii) shows that $\rho_L(F, q) \geq \rho_L(F, p) - 1$, from
which it follows that both inequalities are equalities, and that $y$ is
the minimal element of $F \cap \e^{-1}(q)$. This means that $y \in F_L$, as
desired.

The assertion about Coxeter elements follows from the fact that each
of $F_L$ contains precisely one element with each possible label, which 
completes the proof of (i). The proof of (ii) follows by a symmetrical argument.
\qed\enddemo

\definition{Definition \secd.6}
Let $F$ be a convex subheap of $\heapel$ with full support.
If $F$ has a subheap $F_L$ as in Definition \secd.4, we define the 
{\it lower weight} of $F$ to be the contour of $F_L$.
If $F$ has a subheap $F_H$ as in Definition \secd.4, we define the 
{\it upper weight} of $F$ to be the contour of $F_H$.

If $F$ is the heap of a minuscule element $w$ with full support, then we 
define the lower and upper weights of $w$ to be the lower and upper weights of
$F$ (which will both exist).
\enddefinition

\example{Example \secd.7}
If $F$ is the shaded area in Figure \secba.3, then the lower weight of $F$ is
${+}{-}{-}{+}{+}{-}$, and the upper weight of $F$ is ${+}{-}{+}{-}{+}{+}$.
\endexample

If $\bw$ is a reduced expression a minuscule element $w$ with full support, 
it is possible to define the lower and upper weights of $w$ solely
in terms of the expression $\bw$, as in Remark \secd.3. For example, the 
upper weight of $w$ has a ${+}$ in position $i$ if the leftmost occurrence
of $s_{i-1}$ in $\bw$ occurs to the right of the leftmost occurrence of $s_i$
in $\bw$.

\proclaim{Lemma \secd.8}
Let $\propid$ be the set of proper ideals of
$\heapel$ and let $\Lambda$ be the set of weights.
Let $\psi : \propid \ra \zed \times \Lambda$ be the function that sends a
proper ideal $J \in \propid$ to the pair $(c, \lambda)$, where $(0, 2c)$ is
the maximal element of the vertex chain $J \cap \e^{-1}(0)$, and
$\lambda$ is the upper weight of $J$. Then $\psi$ is a bijection.
\endproclaim

\demo{Proof}
This is a restatement of \cite{{\bf 10}, Proposition 6.4.17 (i)}.
\qed\enddemo

\proclaim{Proposition \secd.9}
The $\asn$-module $V_\propid$ of Proposition \secd.1 is isomorphic
to the module $\asmt$ of Definition \secb.2, where $\asmt$ is regarded
as an $\asn$-module by restriction from $\kqq \otimes_{\kq} \asn$.
An explicit isomorphism is given by 
the $k$-linear map $\theta$ defined by $\theta(b_I) = q^c \lambda$, 
where $b_I  \in V_\propid$, $c \in \zed$
and $\lambda \in \Lambda$ satisfy $\psi(I) = (c, \lambda)$ in the notation
of Lemma \secd.8.
\endproclaim

\demo{Proof}
It is immediate that $\theta$ is a bijection, so it remains to show that
$\theta$ is a module isomorphism.  This is a matter of checking case 
by case that the relations in Proposition \secb.3
are compatible with the bijection $\psi$ of Lemma \secd.8, and the relations
of Proposition \secb.3 are defined in the way that they are precisely so that 
this works.
\qed\enddemo

Recall that each minuscule element of $W$ corresponds to well-defined 
element $u_w \in \asn$, as described following Definition \secba.11.

\proclaim{Lemma \secd.10}
Let $w \in W$ be a minuscule element and let $J, J' \in \propid$. Then the 
following are equivalent:
\item{\rm (i)}{$u_w . b_J = b_{J'}$;}
\item{\rm (ii)}{$J \subseteq J'$ and 
$J' \backslash J$ is isomorphic in $\hcatg$ to the heap $F$ of $w$.}
\endproclaim

\demo{Proof}
The proof is by induction on $r$, where $r = \ell(w)$ is the length of $w$
in (i), and $r$ is the cardinality of $J' \backslash J$ in (ii). If $r = 0$, 
we have $u_w = 1$ and the result follows by taking $J' = J$. If $r = 1$, 
the result follows from Proposition \secd.1. 

We may therefore assume from now on that $r > 1$.
Suppose that $s_{i_1} s_{i_2} \cdots s_{i_r}$ 
is a reduced expression for
$w$, so that $u_w = u_{i_1} u_{i_2} \cdots u_{i_r}$. Let $w' = s_{i_1} w$, so
that $\ell(w') = r-1$. It is immediate from the definitions that $w'$ is
minuscule, and that $u_w = u_{i_1} u_{w'}$.

Suppose that we have $u_w . b_J = b_{J'}$ for some $J, J' \in \propid$.
The above paragraph shows that we have $u_{w'} . b_J \ne 0$, and repeated
applications of Proposition \secd.1 show that $u_{w'} . b_J = b_K$ for some
$K \in \propid$ with $J \subseteq K$. By induction, $K \backslash J$ is the
heap of $w'$. Since $u_w = u_{i_1} u_{w'}$, it follows that 
$u_{i_1} . b_K = b_{J'}$, which means that $J'$ can be obtained from $K$ by
adding a single maximal element $x$ with $\e(x) = i_1$. Since $x \not\in K$,
we also have $x \not\in J$, so $x$ is a maximal element of 
$J' \backslash J = (K \backslash J) \cup \{x\}$. Adding a maximal element
with label $i_1$ to the heap of $w'$ produces the heap of
$w$, proving that (i) implies (ii).

Conversely, suppose that $J \subseteq J'$ and
$J' \backslash J$ is isomorphic in $\hcatg$ to 
the heap $F$ of $w$. Since $w$ has a reduced expression beginning with 
$s_{i_1}$, $F$ has a maximal element $x$ with label $i_1$. Since $J$ and $J'$
are ideals, $x \in J' \backslash J$ is also a maximal element of $J'$. 
It follows that $K := J' \backslash \{x\}$ is a proper ideal that contains
$J$. The construction of $K$ ensures that we have
(a) $K \backslash J = F \backslash \{x\}$, and (b) 
$u_{i_1} . b_K = b_{J'}$.  Using (a), we find that $K \backslash
J$ is the heap of $w'$, which by induction implies that 
$u_{w'} . b_J = b_K$. Applying (b) now completes the proof.
\qed\enddemo

\proclaim{Theorem \secd.11}
\item{\rm (i)}{The set
$\{u_w : w \in W_m\}$ indexed by the minuscule elements of $W$ is a basis
for $\asn$. A word $u_{i_1} u_{i_2} \cdots u_{i_r}$ in the generators is zero
unless $s_{i_1} s_{i_2} \cdots s_{i_r}$ is a reduced expression for a minuscule
element; in particular, any nonzero word in the generating set $u_i$ is equal 
to a basis element.}
\item{\rm (ii)}{The module $V_\propid$ of Proposition \secd.1 is a faithful
$\asn$-module.}
\item{\rm (iii)}{The algebras $\asn$ and $\asn/\langle Q \rangle$ both have
wild representation type.}
\endproclaim

\demo{Proof}
It is clear that the set of all words in the $u_i$ (including the empty word)
form a spanning set for $\asn$. By Proposition \secba.3 (ii), if $\bw$
is not a reduced expression for a minuscule element, then $u_\bw$ is zero.
To prove (i), we therefore need to show that the claimed set is linearly
independent.

By Theorem \secc.13, there is a finite convex subheap $F$ of $\heapel$ 
representing $w$.
By \cite{{\bf 10}, Lemma 3.2.4 (vii)}, a finite convex subheap of a full heap is
the symmetric difference of two nested proper ideals, $J \subseteq J'$, so
by choosing $J$ and $J'$ suitably, we have $F = J' \backslash J$.

Lemma \secd.10 now shows that $u_w . b_J = b_{J'}$. Furthermore, we can recover
$w$ from a knowledge of $J'$ and $J$, because the heap of $w$ is isomorphic to
$J' \backslash J$. It follows that $b_{J'}$ occurs with zero coefficient in
$u_x . b_J$ for every minuscule element $x \ne w$.

Let $u$ be a nonzero element of $\asn$ and choose $w$ so that the coefficient
$c_w$ in $u= \sum_{x \in W_m} c_x u_x$ is nonzero.
The above argument shows that the coefficient of $b_{J'}$ in $u . b_J$
is equal to $c_w$. It follows that the action of $u$ on $V_\propid$ is nonzero.
Since $u$ is an arbitrary nonzero element, this shows that the action is 
faithful, which proves (ii) and also that the set in (i) is linearly 
independent.

To prove (iii), we define the ring $R = k[x, y, z]/I^2$, where
$I$ is the ideal of $k[x, y, z]$ given by $\langle x, y, z \rangle$.
There is a unital homomorphism 
$\pi : \asn \ra R$ of $k$-algebras that sends $u_0$, $u_1$ and $u_2$ to 
$x + I^2$, $y + I^2$, $z + I^2$ respectively, and sends $u_i$ to zero if 
$i > 2$. This homomorphism is well defined because the defining relations 
of $\asn$ all involve words of length $2$ or $3$, and is surjective by
construction. By \cite{{\bf 12}, Theorem 2.10} and \cite{{\bf 11}, Lemma 3}, the algebra
$R$ has wild representation type, as does any algebra having $R$ as a quotient,
which proves that $\asn$ has wild representation type. 

The element $Q$ is a homogeneous polynomial in the $u_i$ of degree $n+1 > 2$,
so we have $\pi(Q) = 0$. The argument of the previous paragraph can then be
applied to $\asn/\langle Q \rangle$, showing that the quotient algebra has
wild representation type, and proving (iii).
\qed\enddemo

\remark{Remark \secd.12}
The argument of Theorem \secd.11 (iii) also applies to any of the generalized
nil Temperley--Lieb algebras of \cite{{\bf 2}}, provided that it has at least 
three generators. In particular, the nil Temperley--Lieb algebras of types
$A$ and affine $A$ have wild representation type, except in trivial cases.
\endremark

\demo{Proof of Theorem \secb.4}
It remains to show that $\asm$ is a faithful $\asn$-module.
Let $u$ be a nonzero element of $\asn$ and choose $w$ so that the coefficient
$c_w$ in $u = \sum_{x \in W_m} c_x u_x$ is nonzero. By the proof of 
Theorem \secd.11, there exist basis elements $b_J$ and $b_{J'}$ of
$V_\propid$ such that $b_{J'}$ appears with nonzero coefficient in $u . b_J$.

Using the isomorphism of Proposition \secd.9, there exist $c, d \in \zed$
and $\lambda, \mu \in \Lambda$ such that the $k$-basis element 
$q^d \mu$ of $\asmt$ appears with nonzero
$k$-coefficient in $u . (q^c \lambda)$. If $e \in \zed$ is taken large enough
that $c + e$ and $d + e$ are both nonnegative, then $q^{c+e} \lambda$ and
$q^{d+e} \mu$ are $k$-basis elements of $\asm$ with the property that 
$q^{d+e} \mu$
appears with nonzero coefficient in $u . (q^{c+e} \lambda)$. It follows that
the actions of $\asn$ on $\asm$ and on $\asmt$ are faithful, which completes the
proof.
\qed\enddemo

\head \sece. $\asn$ as a ring of matrices \endhead

The main purpose of Section \sece\ is to prove Theorem \secb.6. 
The key to this is a certain central element $Q \in \asn$, which we 
define as follows. 

\definition{Definition \sece.1}
If $\lambda \in \Lambda$ is a weight, and $w(\lambda)$ is Coxeter element
corresponding to $\lambda$ (as in Remark \secd.3), we write $u_\lambda$
for $u_{w(\lambda)}$. We also define $Q \in \asn$ by $$
Q = \sum_{\lambda \in \Lambda} u_\lambda
.$$
\enddefinition

It will turn out that $Q$ is closely related to the following automorphism
of $\heapel$.

\proclaim{Lemma \sece.2}
The map $\tau : \heapel \ra \heapel$ given by $\tau((c, d)) = ((c, d+2))$ is
an automorphism in $\hcatg$.
\endproclaim

\demo{Proof}
This is the automorphism described in \cite{{\bf 10}, Lemma 5.3.1}. (It is also not
hard to check the details directly.)
\qed\enddemo

\proclaim{Lemma \sece.3}
Let $w \in W$ be minuscule and let $\lambda \in \Lambda$.
\item{\rm (i)}{If $u_\lambda u_w$ is nonzero, then there exists a unique
$\mu \in \Lambda$ for which $u_\lambda u_w = u_w u_\mu$.}
\item{\rm (ii)}{If $u_w u_\lambda$ is nonzero, then there exists a unique
$\mu \in \Lambda$ for which $u_w u_\lambda = u_\mu u_w$.}
\endproclaim

\demo{Proof}
We first prove (i). Suppose that $u_\lambda u_w$ is nonzero, and let $w(\lambda)$
be the Coxeter element corresponding to $\lambda$. By Theorem \secd.11 (i),
$x = w(\lambda) w$ is a minuscule element, and $\ell(x) = \ell(w(\lambda))
+ \ell(w)$. Since $w(\lambda)$ has full support, so does $x$. 
By Lemma \secba.13, the heap of $x$, which we identify with a convex subset
of $\heapel$ by Theorem \secc.13, has a filter isomorphic to the heap of
$w(\lambda)$, which means that the upper weight of $x$ is $\lambda$.

If it is possible to express $u_\lambda u_w$ in the form $u_y u_\mu$ for some
$y \in W_m$, an argument like that of the previous paragraph shows that 
$\mu$ must be the lower weight of $x$, which proves the assertion about
uniqueness. Lemma \secba.13 shows that
the heap of $x$ has an ideal isomorphic to the heap of $w(\mu)$, and a
filter isomorphic to the heap of $y$. 

It remains to show that $w$ and $y$ have isomorphic heaps. Suppose that
$F$ is the heap of $x$ and that $F \cap \e^{-1}(\{s, t\})$ is an edge chain
in $F$. By construction, the corresponding $\{s, t\}$-edge chain of the heap
of $w$ (respectively, the heap of $y$) is obtained by removing the two highest
(respectively, two lowest) elements of the edge chain of $F$. By Lemma
\secc.6, all three edge chains are alternating. It follows that the 
$\{s, t\}$-edge chain of $y$ (which may be the empty set) can be obtained 
from that of $w$ by applying the automorphism $\tau$ of Lemma \sece.2.
Since each edge chain behaves in the same way, the entire heap of $y$
can be obtained from that of $w$ by applying $\tau$.
This shows that the heaps are
isomorphic, and therefore that $y = w$, as required.

The proof of (ii) follows by a symmetric argument, by reversing
the order of multiplication.
\qed\enddemo

\proclaim{Corollary \sece.4}
The element $Q \in \asn$ lies in the centre of $\asn$.
\endproclaim

\demo{Proof}
Lemma \sece.3 shows that $Q$ commutes with $u_w$ for all $w \in W_m$, and
these form a basis for $\asn$ by Theorem \secd.11 (i).
\qed\enddemo

\proclaim{Lemma \sece.5}
Let $\lambda \in \Lambda$ and let $\tau$ be as in Lemma \sece.2. 
The action of $u_\lambda$ on $\propid$ is given by $$
u_\lambda . b_J = \cases
\tau(b_J) & \quad \text{ if $\lambda$ is the upper weight of $J$}, \cr
0 & \quad \text{ otherwise.}
\endcases
$$
\endproclaim

\demo{Proof}
Since $u_\lambda$ is a word in the generators, $u_\lambda . b_J$ must either be
zero, or equal to $b_K$ for some $K \in \propid$ that contains $J$.
Suppose that $u_\lambda . b_J = b_K$. By Lemma \secd.10, $K \backslash J$ is
isomorphic to the heap of $u_\lambda$, which means that $K \backslash J$
contains precisely one element of each label. Since $K$ is an ideal, this means
that we can obtain $K$ from $J$ by adjoining, for each label $p$, the minimal
element of the chain $(\heapel \backslash J) \cap \e^{-1}(p)$. The automorphism
$\tau$ has exactly the same effect on $J$, which forces $K = \tau(J)$.

Let $\mu$ be the upper weight of $J$, which exists because $J$ is a 
proper ideal. It follows that $J$ has a filter isomorphic to the heap $F$ of
$u_\mu$; note that $F$ is also isomorphic to $\tau(J) \backslash J$.
Combined with the previous paragraph, this means that $\mu = \lambda$ and
that $u_\mu . b_J = b_{\tau(J)}$, which completes the proof.
\qed\enddemo

\proclaim{Lemma \sece.6}
Let $V_\propid$ be the $\asn$-module of Proposition \secd.1.
\item{\rm (i)}{The element $Q$ acts on $V_\propid$ via $Q . b_J = b_{\tau(J)}$.}
\item{\rm (ii)}{The element $Q$ acts on $\asmt$, and on $\asm$, as
multiplication by $q$.}
\endproclaim

\demo{Proof}
The element $Q$ is the sum of the $u_\lambda$ for all possible weights 
$\lambda$. Lemma \sece.5 shows that exactly one of these terms will act in
a nonzero way on $b_J$, and the one that does not act as zero will send
$b_J$ to $b_{\tau(J)}$, proving (i).

The identification between $\asmt$ and $V_\propid$ of Proposition \secd.9
shows that the automorphism $\tau$ induces the map of multiplication by $q$
on $\asmt$, and on $\asm$. Part (ii) follows from this.
\qed\enddemo

\proclaim{Lemma \sece.7}
For $\lambda, \mu \in \Lambda$, we have $$
u_\lambda u_\mu = \delta_{\lambda, \mu} Q u_\mu
$$ in $\asn$, where $\delta$ is the Kronecker delta.
\endproclaim

\demo{Proof}
Since $\asn$ acts faithfully on $V_\propid$, it is enough to show that the
two sides of the equation act in the same way on a given $b_J \in V_{\propid}$.

If the upper weight of $b_J$ is not $\mu$, then $u_\mu . b_J = 0$ by 
Lemma \sece.5 and both
sides act as zero, which satisfies the statement.

If the upper weight of $b_J$ is $\mu$, then $u_\mu . b_J = b_{\tau(J)}$ by
Lemma \sece.5.
Since the isomorphism $\tau$ preserves upper weights, $\tau(J)$ also has
an upper weight of $\mu$. If $\lambda \ne \mu$, Lemma \sece.5 shows that 
$u_\lambda . b_{\tau(J)} = 0$, and both sides of the equation act as zero, 
again satisfying the statement.

Finally, if $\lambda = \mu$, iterating the argument of the previous paragraph
shows that $u_\lambda u_\mu . b_J = b_{\tau^2(J)}$, and Lemma \sece.6 (i) shows
that $Q u_\mu . b_J = b_{\tau^2(J)}$, which completes the proof in this case.
\qed\enddemo

\proclaim{Lemma \sece.8}
Let $\rho$ be the usual rank function on $\heapel$ and 
let $I_1$ and $I_2$ be (not necessarily finite) convex subheaps of $\heapel$
with full support, such that each vertex chain of $I_1$ (respectively, $I_2$)
is bounded above (respectively, below). Suppose that the upper weight of
$I_1$ is equal to the lower weight of $I_2$.  Then there exists an integer
$r$ such that 
$F = I_1 \cup \tau^r(I_2)$
is a convex subheap of $\heapel$ having $I_1$ as an ideal and $\tau^r(I_2)$ 
as a filter.
\endproclaim

\demo{Proof}
Denote the common weight mentioned the statement by $\lambda$. By applying
a suitable power of $\tau$ to $I_2$, we may ensure that 
$\rho_L(I_2, p) = \rho_H(I_1, p) + 2$ for each label $p$. We choose $r \in \zed$
so that this is the case.

Since $F$ has full support, it suffices
by Theorem \secc.1 to show that for each edge $\{s, t\}$ of
$\Gamma$, the edge chain $I \cap \e^{-1}(\{s, t\})$ is a convex subset of
$E \cap \e^{-1}(\{s, t\})$. The corresponding property holds for $I_1$ and 
$\tau^r(I_2)$
individually by Remark \secc.2, so we will be done if we can show that 
the lowest element, $x$, of the edge chain
$\tau^r(I_2) \cap \e^{-1}(\{s, t\})$
covers
the highest element, $y$, of the edge chain 
$I_1 \cap \e^{-1}(\{s, t\})$.

Without loss of generality, we can take $s = p$ and $t = p+1$. There are two
cases to consider. If $\lambda$ has a ${+}$ in the $p$-position, then
we have $\e(y) = p+1$ and $\e(x) = p$. Furthermore, the highest element
$z$ in $I_1$ such that $\e(z) = p$ satisfies $z \prec y$ and 
$\e(z) = p$. By construction, we have $\rho(x) = \rho(z) + 2$, and 
$z \prec y$ implies that $\rho(y) = \rho(z) + 1$. It follows that we
have $y = (p+1, c)$ and $x = (p, c+1)$ for some $c \in \zed$, meaning that
$x$ covers $y$, as required.

The other case, in which $\lambda$ has a ${-}$ in the $p$-position, is dealt
with by a symmetric argument.
\qed\enddemo

\definition{Definition \sece.9}
If $w$ is a minuscule element with full support, such that $\lambda$ is the
upper weight of $w$, $\mu$ is the lower weight of $w$, and there are $r$
occurrences of the generator $s_0$ in $w$, then we write
$C_{\lambda, \mu}^r$ for $w$.
\enddefinition

\example{Example \sece.10}
The convex subheap of Figure \secba.3 has one element with label $0$, so 
it corresponds to the element $C_{\lambda, \mu}^1$, where
$\mu = {+}{-}{-}{+}{+}{-}$, and $\lambda = {+}{-}{+}{-}{+}{+}$.
\endexample

Not every choice of $\lambda$, $\mu$ and $r \geq 0$ in Definition \sece.9 
corresponds to a minuscule element. However, the notation is unique when it
does apply, as the following result shows.

\proclaim{Lemma \sece.11}
Maintain the notation of Definition \sece.9.
\item{\rm (i)}{If $C_{\lambda, \mu}^r$
corresponds to a minuscule element with full support, then
there is a unique minuscule element with this property.}
\item{\rm (ii)}{For every $\lambda, \mu \in \Lambda$, there exists 
$r \geq 1$ such that $C_{\lambda, \mu}^r$ corresponds to a minuscule element.}
\endproclaim

\demo{Proof}
We will show how to recover the heap $F$ of $w$ up to isomorphism, considered
as a convex subheap of $\heapel$, as in Theorem \secc.13. By applying a power
of $\tau$ if necessary, we may assume that the lowest element of the vertex
chain $F \cap \e^{-1}(0)$ has rank $1$. 
The rank of the highest element of $F \cap \e^{-1}(0)$ is then $r$.
The other minimal elements of the vertex
chains of $F$ are then determined by $\mu$, and the other maximal 
elements of $F$ are determined by $\lambda$. The rest of the heap can be
filled in because each vertex chain of $F$ is an interval in the 
corresponding vertex chain of $\heapel$, and this establishes (i).

To prove (ii), let $F_1$ and $F_2$ be convex subheaps of $\heapel$ 
representing $u_\lambda$ and $u_\mu$, respectively. By applying a sufficiently
large power of $\tau$ to $F_1$, we may assume that every element of $F_1$
has greater rank than every element of $F_2$. We may now fill in the vertex
chains as in the construction used to prove (i) above to produce a finite
convex subheap with lower weight $\mu$ and upper weight $\lambda$. This
corresponds to a minuscule element with the required properties.
\qed\enddemo

Lemma \sece.11 ensures that the statement of the next result is well defined.

\proclaim{Proposition \sece.12}
\item{\rm (i)}{If $C_{\lambda, \mu}^r$ and $C_{\nu, \xi}^s$ are 
minuscule elements with full
support, then we have $r, s \geq 1$ and $$
C_{\lambda, \mu}^r C_{\nu, \xi}^s = \delta_{\mu, \nu} C_{\lambda, \xi}^{r+s}
,$$ where $\delta$ is the Kronecker delta.}
\item{\rm (ii)}{If $C_{\lambda, \mu}^r$ corresponds to a minuscule element 
with full support, then we have $Q C_{\lambda, \mu}^r = C_{\lambda, \mu}^{r+1}$.}
\endproclaim

\demo{Proof}
We first prove (i).
Any minuscule element with full support must have the label $0$ in its support,
and this shows that $r, s \geq 1$.
The hypotheses about weights show that we can write $C_{\lambda, \mu} = 
u_w u_\mu$ and $C_{\nu, \xi} = u_\nu u_x$ for suitable minuscule elements $w$
and $x$. If $\mu \ne \nu$ then Lemma \sece.7 shows that the product in the 
statement is zero, as required.

For the rest of the proof of (i), we may assume that $\mu = \nu$. Let $F_1$ 
(respectively, $F_2$)
be a convex subheap of $\heapel$ representing $C_{\nu, \xi}$ (respectively,
$C_{\lambda, \mu}$). The hypotheses show that the upper weight of $F_1$ is 
equal to the lower weight of $F_2$. By applying a suitable power of $\tau$
to $F_1$, we may apply Lemma \sece.8 with $I_1 = F_1$ and $I_2 = F_2$ to
obtain a finite convex subheap $I$ with $I_1$ as an ideal and $I_2$ as a filter.

The lower weight of $I$ is the same as that of $I_1$, the upper weight of
$I$ is the same as that of $I_2$, and the number of elements labelled $0$ in
$I$ is equal to $r+s$. The proof of (i) is now completed by Lemma \secba.13.

The element $C_{\lambda, \mu}^r$ in the statement of (ii) 
has upper weight $\lambda$, so we have 
$C_{\lambda, \mu}^r = u_\lambda u_x$ for some minuscule $x$. By Lemma \sece.7,
we then have $$
Q C_{\lambda, \mu}^r = Q u_\lambda u_x = u_\lambda^2 u_x = 
u_\lambda C_{\lambda, \mu}^r 
.$$ Since $u_\lambda = C_{\lambda, \lambda}^1$, it follows from part (i) that
$u_\lambda C_{\lambda, \mu}^r = C_{\lambda, \mu}^{r+1}$, which proves (ii).
\qed\enddemo

\proclaim{Lemma \sece.13}
Let $\mlqq$ be the full matrix ring whose rows and columns are
indexed by the set of weights $\Lambda$, and whose entries are Laurent 
polynomials in $q$.
\item{\rm (i)}{The 
set of all $u_w$, as $w$ ranges over all minuscule elements with full
support, forms an ideal, $I$, of $\asn$.}
\item{\rm (ii)}{There is an injective homomorphism $\kappa_0$
of $k$-algebras from $I$ to $\mlqq$ 
satisfying $\kappa_0(C_{\lambda, \mu}^r) = q^r E_{\lambda, \mu}$,
where $E_{\lambda, \mu}$ is the matrix unit with $1$ in the $(\lambda, \mu)$
position and zeros elsewhere, and $r \geq 1$.}
\item{\rm (iii)}{The map $\kappa_0$ sends $u_\lambda$ to $q E_{\lambda, \lambda}$,
and sends $Q$ to $q I_\Lambda$, where $I_\Lambda$
is the identity matrix.}
\item{\rm (iv)}{The homomorphism of (ii) extends to an injective, unital 
homomorphism $\kappa$ of $k$-algebras from $\asn$ to $\mlqq$ defined by 
$\kappa(u) = q^{-1} \kappa_0(Qu)$.}
\item{\rm (v)}{The image of $\kappa$ lies in the subring $\mlq$.}
\endproclaim

\demo{Proof}
Part (i) follows from Theorem \secd.11 (i), and part (ii) is a restatement
of Proposition \sece.12 (i).

Part (iii) follows from (ii), combined with the definition of $Q$ and the
fact that $u_\lambda = C_{\lambda, \lambda}^1$.

Since each $u_\lambda$ lies in $I$, it follows that $Q$ and $Qu$ also lie in
$I$, so $\kappa(u)$ is well defined.  By (ii) and 
Corollary \sece.4, we now have $$\kappa_0(Qu) \kappa_0(Qv) = \kappa_0(Q^2uv).$$
Because $Quv$ lies in $I$, Proposition \sece.12 (ii) shows that 
$\kappa_0(Q^2uv) = q \kappa_0(Quv)$. Combining these equations, we then have
$(q^{-1} \kappa_0(Qu))(q^{-1} \kappa_0(Qv)) = q^{-1} \kappa_0(Quv)$, which
shows that $\kappa$ is a homomorphism. Finally, part (ii) shows that $\kappa$
is unital, which proves (iv).

Part (v) follows from (iv) and the fact that the integers $r$ in (ii) satisfy
$r \geq 1$.
\qed\enddemo

\definition{Definition \sece.14}
The $k$-algebra $\asnt$ is defined by adjoining a new generator, $Q^{-1}$,
to the presentation in Definition \secb.1, together with the new relations $$
Q Q^{-1} = Q^{-1} Q = 1
,$$ where $Q$ is the element of $\asn$ defined in Definition \sece.1.
\enddefinition

\proclaim{Lemma \sece.15}
The natural homomorphism from $\asn$ to $\asnt$ is injective.
\endproclaim

\demo{Proof}
By Lemma \sece.6 (i), $Q$ acts invertibly on $V_\propid$. It
follows from the presentation in Definition \sece.14 that the representation
$\sigma : \asn \ra \End_k(V_\propid)$ can be extended 
to a representation $\sigma' : \asnt \ra \End_k(V_\propid)$ by defining $$
Q^{-1} . b_J = b_{\tau^{-1}(J)}
$$ for all $J \in \propid$. By Theorem \secd.11 (ii), $\sigma$ is faithful,
which means that
$\sigma'$ is faithful when restricted to $\asn$, and this can only happen if
the map in the statement is injective.
\qed\enddemo

\demo{Proof of Theorem \secb.6}
We will prove that the injective, unital homomorphism 
$\kappa : \asn \ra \mlqq$ in Lemma
\sece.13 (iv) extends to an isomorphism, $\kappa'$, of $k$-algebras from 
$\asnt$ to $\mlqq$.

By Lemma \sece.13 (iii), the element $Q \in \asn$ is mapped to the unit
$q I_\Lambda$ of $\mlqq$, so we can extend $\kappa$ to a homomorphism
$\kappa' : \asnt \ra \mlqq$ by defining $\kappa'(Q^{-1}) = q^{-1} I_\Lambda$.

Suppose that $v \in \ker(\kappa')$.
Because $\asnt$ is a localization of $\asn$, it follows that we have
$Q^r v \in \asn$ for some positive integer $r$. Since $Q^r v$ lies in the
intersection of $\asn$ with the ideal $\ker(\kappa')$, the injectivity of
$\kappa$ proves that $Q^r v = 0$, and 
therefore that $v = 0$. This proves that $\kappa'$ is injective.

Let $\lambda, \mu \in \Lambda$. By Lemma \sece.11 (ii), there exists
an integer $r$ such that $\kappa(u) = q^r E_{\lambda, \mu}$ for some 
$u \in \asn$. Given any $s \in \zed$, it follows that 
$\kappa'(Q^{s-r}u) = q^s E_{\lambda, \mu}$. 
Since the set $$
\{ q^s E_{\lambda, \mu} : s \in \zed, \lambda \in \Lambda, \mu \in \Lambda \}
$$ is a $k$-basis for $\mlqq$, it follows that $\kappa'$ is surjective.
\qed\enddemo

\head \secf. Representation theory \endhead

In Section \secf, we will prove Theorem \secb.5, concerning the centre
of $\asn$, and Theorem \secb.8, describing the representation theory of
$\asn$.
The construction in Section \sece\ of $\asn$ in terms of matrices makes many
of the structural features of $\asn$ more transparent. In particular,
we can describe the centre of $\asn$, as follows.

\proclaim{Proposition \secf.1}
The centre of $\asn$ is equal to $k[Q]$.
\endproclaim

\demo{Proof}
It follows from Corollary \sece.4 that $k[Q]$ is contained in the centre of
$\asn$, so it remains to show that the centre is no larger than this.

Suppose that $u$ is an element of the centre of $\asn$. Since $u$ commutes with
$Q$, it also commutes with $Q^{-1}$ when regarded as an element of $\asnt$,
which means that $u$ lies in the centre of $\asnt$. By Theorem \secb.6,
$\asnt$ is isomorphic to $\mlqq$, and the centre of a full matrix ring over
a commutative ring consists precisely of the scalar multiples of the identity
matrix. By Lemma \sece.13 (v), this scalar multiple must lie in $k[q]$, 
which shows that $u \in k[Q]$, as required.
\qed\enddemo

\proclaim{Lemma \secf.2}
The $k$-algebras $\kqq \otimes_\kq \asn$ and $\asnt$ are isomorphic, where
we regard $\asn$ as a left $\kq$-module via $q.u = Qu$.
\endproclaim

\demo{Proof}
The map $\kappa$ of Lemma \sece.13 (v) identifies $\asn$ with a $k$-subalgebra
$A$ of $\mlq$, and Lemma \sece.13 (iii) shows that this identification 
respects the module structure in the statement. 

It follows by extension of scalars that we have 
$\mlqq \cong \kqq \otimes_\kq \mlq$. 
Since $\kqq$ is a localization of $\kq$, it follows that $\kqq$ is flat
as a $\kq$-module, which shows that $\kqq \otimes_\kq A \leq \kqq \otimes_\kq
\mlq$.
Lemma \sece.11 (ii) shows that for
any $\lambda, \mu \in \Lambda$, $A$ contains an element of the form
$q^r E_{\lambda, \mu}$. This means that we have 
$\kqq \otimes_\kq A = \kqq \otimes_\kq \mlq$, from which the statement follows.
\qed\enddemo

\demo{Proof of Theorem \secb.5}
The assertions about $Q$ and the centre were proved in Lemma \sece.13 (iii)
and Proposition \secf.1. It remains to show that $\asn$ is a free $k[Q]$-module
of rank $2^{2n} = |\Lambda|^2$.

By Lemma \sece.13 (v), $\asn$ is isomorphic to a $k$-subalgebra $A$ of $\mlq$.
It follows that $A$ is a $\kq$-submodule of the finitely generated 
$\kq$-module $\mlq$, and because $\kq$ is Noetherian, $A$ is finitely generated
as a $\kq$-module.

Because $\mlqq$ is a finitely generated free $\kqq$-module and $\kqq$ is a
domain, $\mlqq$ is torsion free as a $\kqq$-module. By restriction, it follows
that $A$ is a torsion free $\kq$-module. Since $\kq$ is a principal
ideal domain and $A$ is finitely generated as a $\kq$-module, it follows that 
$A$ is a free $\kq$-module.

Since $\kqq$ is flat as a
$\kq$-module, the rank of $A$ as a free $\kq$-module is equal
to the rank of $\mlqq$ as a free $\kqq$-module. The latter rank is 
equal to $|\Lambda|^2 = 2^{2n}$, and this completes the proof.
\qed\enddemo

\proclaim{Lemma \secf.3}
The set $$
\{e_\lambda : \lambda \in \Lambda\}
,$$ where $e_\lambda := Q^{-1} u_\lambda$, is a decomposition 
of the identity element of $\asnt$ into orthogonal idempotents. Each element
$e_\lambda$ generates $\asnt$ as a two-sided ideal.
\endproclaim

\demo{Proof}
By Theorem \secb.6, $\asnt$ is isomorphic to $\mlqq$. By Lemma \sece.13 (iii),
this isomorphism identifies the elements $e_\lambda$ with the diagonal 
matrix units $E_{\lambda, \lambda}$ of $\mlqq$. The assertions now follow by
standard properties of matrix rings.
\qed\enddemo

We now recall some well known properties of Morita equivalence. Suppose that
$A$ is an algebra and $e$ is an idempotent of $A$ such that $A = AeA$.
In this situation, $A$ and the algebra $eAe$ have equivalent
categories of left modules, and an equivalence can be induced by the functors
$G_1 : A\text{-}\Mod \ra eAe\text{-}\Mod$ and 
$G_2 : eAe\text{-}\Mod \ra A\text{-}\Mod$ 
defined by $G_1(M) = eM$ and $G_2(N) = Ae \otimes_{eAe} N$.

\proclaim{Lemma \secf.4}
Fix $\lambda \in \Lambda$.
The $\asnt$-module $\asmt$ is isomorphic, as a left $\asnt$-module, to the 
left ideal $\asnt e_\lambda$.
\endproclaim

\demo{Proof}
By Lemma \secf.3, we may use Morita equivalence to reduce the problem to proving
that $e_\lambda \asmt$ and $e_\lambda \asnt e_\lambda$ are isomorphic
as $B$-modules, where $B = e_\lambda \asnt e_\lambda$. 

The set $e_\lambda \asmt$
consists of all $\kqq$-multiples of the weight $\lambda$. The isomorphism
of Lemma \sece.13 identifies $B$ with the $\kqq$-multiples of the matrix
unit $E_{\lambda, \lambda}$. Both modules are isomorphic as $B$-modules to $B$
itself, and this completes the proof.
\qed\enddemo

The next definition is reminiscent of that of the $\asnt$-modules defined
in Section \secb.

\definition{Definition \secf.5}
If $c \in k \backslash \{0\}$ and $m \in \enn$, we define 
the $\asnt$-module 
$\asmlt{c}{m}$ to be $$
\asmlt{c}{m} := \asmt \otimes_{\kqq} 
\frac{\kqq}{\langle (q - c)^m \rangle}
,$$ where the right action of $Q$ on $\asmt$ is by multiplication by $q$.
\enddefinition

\proclaim{Lemma \secf.6}
Suppose that the field $k$ is algebraically closed, and let $c \in k
\backslash \{0\}$ and $m \in \enn$.
The modules $\asmlt{c}{m}$ are indecomposable, and every finite dimensional
indecomposable $\asnt$-module is isomorphic to one of the modules 
$\asmlt{c}{m}$, for a unique value of $c$ and $m$. The module $\asmlt{c}{m}$
is irreducible if and only if $m = 1$.
\endproclaim

\demo{Proof}
Fix a weight $\lambda \in \Lambda$.
By Lemma \secf.4, $\asmt$ is isomorphic, as a left module, to $\asmt e_\lambda$.
Because $Q$ is central and $e_\lambda$ is idempotent, the right action of
$\kqq$ on $\asmt$ is compatible with the right action of 
$e_\lambda k[Q, Q^{-1}] e_\lambda = k[Q, Q^{-1}] e_\lambda$ 
on $\asmt e_\lambda$, provided that we 
identify $\kqq$ and $e_\lambda k[Q, Q^{-1}] e_\lambda$ in the obvious way.

With these identifications, the module $\asmlt{\l}{m}$ corresponds under
the Morita equivalence to the $\kqq$-module $$
\frac{\kqq}{\langle (q - c)^m \rangle}
,$$ which is indecomposable in general and irreducible if and only if $m = 1$.
Because $\kqq$ is a principal ideal domain, every finite dimensional
indecomposable $\kqq$-module is of the form $\kqq/\langle f(x) \rangle$, where
$f(x) \in \kqq$ is a power of an irreducible polynomial $p(x)$. Because $k$ is 
algebraically closed, $p(x)$ must be a unit multiple of $x - c$ for
some $c \in k \backslash 0$. Since Morita equivalence respects 
indecomposability and irreducibility, the result follows.
\qed\enddemo

\proclaim{Lemma \secf.7}
Let $M$ be a finite dimensional vector space over $k$.
\item{\rm (i)}{If $\sigma : \asnt \ra \End_k(M)$ is a representation of
$\asnt$, then the restriction of $\sigma$ to $\asn$ is a representation of
$\asn$ with the same image.}
\item{\rm (ii)}{If $\sigma : \asn \ra \End_k(M)$ is a representation of
$\asn$ in which $Q \in \asn$ acts invertibly, then $\sigma$ can be extended
to a representation of $\asnt$ with the same image.}
\endproclaim

\demo{Proof}
Suppose that $\sigma$ satisfies the conditions of (i). Since $Q \in \asnt$ is
a unit, the minimal polynomial of the action of $Q$ on $M$ has nonzero
constant term. It follows that $\sigma(Q^{-1})$ is a polynomial in $\sigma(Q)$.
This shows that the restriction of $\sigma$ to $\asn$ has the same image
as $\sigma$, and proves (i).

Now suppose that $\sigma$ satisfies the conditions of (ii). Since $Q$ acts
invertibly, $\sigma$ can be extended to a representation of $\asnt$
on the same module. Applying (i) to this extended representation shows
that the representations have the same image, which proves (ii).
\qed\enddemo

\proclaim{Lemma \secf.8}
Fix $r \in \enn$.
\item{\rm (i)}{All but finitely many of the basis elements
$\{u_w : w \in W_m\}$ lie in the ideal $\langle Q^r \rangle$ of $\asn$.}
\item{\rm (ii)}{The images of the generators $u_i$ of $\asn$ generate
a nilpotent (two-sided) ideal of codimension $1$ in the finite dimensional
$k$-algebra $\asn/\langle Q^r \rangle$.}
\endproclaim

\demo{Proof}
If $w \in W_m$ has full support and upper weight $\lambda$, then we can
write $u_w = u_\lambda u_{w_1}$. If $w_1$ also has full support, Lemma \sece.7
shows that the upper weight of $w_1$ is also equal to $\lambda$. Iterating this,
we can uniquely factorize $u_w = u_\lambda^c u_x$, where $c \in \enn$ and 
$x \in W_m$ does not have full support. If $c > 1$, we may apply Lemma \sece.7
again to write $u_w = Q^{c-1} u_\lambda u_x$.

Suppose that $u_w$ does not lie in the ideal generated by $Q^r$.
We need to have $c \leq r$, so there are finitely many choices for $c$.
Since $W$ is an affine Weyl group, any subgroup of $W$ generated by a proper
subset of the generators $S$ is finite. It follows that $W$ has finitely many
elements without full support, minuscule or otherwise, so there are finitely
many choices for $x$. Finally, $\Lambda$ is finite, so there are finitely
many choices for $\lambda$. This shows that there are finitely many choices
for $w$, proving (i).

By part (i) and Theorem \secd.11 (i), any sufficiently long word in
the generators $u_i$ lies in the ideal generated by $Q^r$, and this implies
that the quotient is finite dimensional. It follows that
the image of the codimension 1 ideal of $\asn$ generated by the $u_i$ maps 
to a nilpotent ideal in the quotient, and (ii) follows.
\qed\enddemo

\proclaim{Lemma \secf.9}
Let $f(q) \in \kq$ be a polynomial with nonzero constant term. Then the natural
map $\gamma: \kq/\langle f(q) \rangle \ra \kqq/\langle f(q) \rangle$ is an 
isomorphism of $k$-algebras.
\endproclaim

\demo{Proof}
The ideal generated by $f(q)$ in $\kq$ is contained in the ideal generated
by $f(q)$ in $\kqq$, so $\gamma$ is well-defined.

Let $g(q) \in \kqq$. Since $f(q)$ has nonzero constant term, we can only
have $g(q)f(q) \in \kq$ if $g(q) \in \kq$. It follows from this that 
$\ker \gamma = 0$ and $\gamma$ is injective.

The fact that $f(q)$ has nonzero constant term also shows that 
$q + \langle f(q) \rangle$ is a unit in $\kq/\langle f(q) \rangle$, which
implies that $\gamma$ is surjective.
\qed\enddemo

\proclaim{Lemma \secf.10}
For any $c \in k \backslash \{0\}$ and let $m \in \enn$, 
the $\asnt$-module $\asmlt{c}{m}$, when regarded as a $\asn$-module by
restriction, is isomorphic to the module $\asml{c}{m}$ of Section \secb.
\endproclaim

\demo{Proof}
The $\asn$-actions on the two modules both follow the formula in 
Proposition \secb.3, except that 
in $\asmlt{c}{m}$, the coefficients involving $q$
lie in $\kqq/\langle (q-c)^m \rangle$, and
in $\asml{c}{m}$, the coefficients involving $q$
lie in $\kq/\langle (q-c)^m \rangle$. 
Lemma \secf.9 shows that these rings are canonically isomorphic, from which
the result follows.
\qed\enddemo

\demo{Proof of Theorem \secb.8}
The assertion about dimension in part (i) 
of the theorem follows from the definition of $\asml{c}{m}$, together
with the facts that $|\Lambda| = 2^n$ and $\kq/\langle (q-c)^m \rangle$ 
has dimension $m$.

Lemma \secf.7 gives a canonical identification between finite dimensional
$\asn$-modules on which $Q$ acts invertibly, and finite dimensional
$\asnt$-modules. Because this correspondence preserves the image of the
representation, it respects indecomposability and irreducibility.

Let $M$ be a finite dimensional indecomposable $\asn$-module over $k$ on
which $Q$ acts invertibly. 
By Lemma
\secf.6, $M$ is isomorphic to a unique module of the form $\asmlt{c}{m}$,
and this module is irreducible if and only if $m = 1$.
Lemma \secf.10 then shows that $M$ is isomorphic to a unique module of the form
$\asml{c}{m}$, which again is irreducible if and only if $m = 1$.

Suppose from now on that $Q$ does not act invertibly. Since $M$ is finite
dimensional and $k$ is algebraically closed, the Jordan canonical form of
the action of $Q$ on $M$ shows that $Q$ acts nilpotently. The action of
$\asn$ on $M$ then factors through a quotient algebra of the form 
$A = \asn/\langle Q^r \rangle$. Lemma \secf.8 (ii) shows that the
images of the generators $u_i$ generate the Jacobson radical of $A$,
which is nilpotent and has codimension $1$. This shows
that the only irreducible module for $A$ is the one dimensional module
on which all generators $u_i$ act as zero.
\qed\enddemo

\head Concluding remarks \endhead

As we have mentioned, it is possible to define an algebra of raising operators
using any full heap (in the sense of \cite{{\bf 10}}). However, we have concentrated
on type $\twd$ in this paper because Dynkin diagrams of types other than
$\twd$ and $A_n^{(1)}$ do not behave as well. A large part of the reason
for this is that in these two types, the lowest positive imaginary root 
associated to the affine Kac--Moody algebra contains each fundamental root
with coefficient $1$, but this is not true in other affine types. 

The fact that all the coefficients above are equal to $1$ creates interesting
connections with Coxeter elements in type $\twd$.`<
For example, Theorem \secc.13 can be
reformulated as saying that for any reduced expression $\bx$ for a fixed 
Coxeter  element $x$, and every minuscule element $w$, there is a reduced 
expression for $x^r$ (for some integer $r$) that contains as a subword a 
reduced expression for $w$. Furthermore, all subwords of $\bx^r$ are reduced
expressions for minuscule elements.

Another positive feature shared by 
$\asn$ and the nil Temperley--Lieb algebra in type
affine $A$ is that the particle configuration representations are rich enough
that the words $u_w$, with $w$ minuscule, act linearly independently. 
It turns out that this does not happen for representations arising from other
full heaps. More specifically, for each such representation, there is a product
of distinct commuting generators $u_i$ in the algebra that acts as zero.

Many of the techniques of this paper can be adapted to work for the 
affine nil Temperley--Lieb algebra of type $A$. In particular, the
``fermionic representations'' approach in \cite{{\bf 13}} can be constructed from 
scratch using heaps. Furthermore, adding additional central elements to
this algebra will result in a direct sum of matrix rings over rings of Laurent
polynomials, and families of indecomposable representations can be constructed
from this using Morita equivalence. It is not possible to do much better than
this, because the type affine $A$ algebras also have wild representation type.

\head Acknowledgements \endhead

I thank Sarah Salmon and the anonymous referee for suggesting improvements and
pointing out some errors in earlier versions of this paper.

\leftheadtext{} \rightheadtext{}
\vfill\eject

\end